%% file: torsion.tex
\documentclass[11pt]{amsart}

\usepackage{amsmath, amsthm, amssymb}

\usepackage{mathtools}
\usepackage{leftindex}
\usepackage{subcaption}
\usepackage{bm}

\usepackage[T1]{fontenc}                    
\usepackage{lmodern}
\usepackage{palatino}                       
\usepackage{inconsolata}                    
\usepackage{comment}

\usepackage[abs]{overpic}		  
\usepackage{import}
\usepackage{transparent}

\usepackage{xcolor}
\definecolor{lred}{RGB}{226, 106, 106}
\definecolor{nred}{RGB}{237, 28, 36}
\definecolor{ddred}{RGB}{255, 0, 0}
\definecolor{lblue}{RGB}{52, 152, 219}
\definecolor{nblue}{RGB}{0, 174, 239}
\definecolor{lyellow}{RGB}{232, 197, 91}
\definecolor{dgreen}{RGB}{0, 148, 68}
\definecolor{l1yellow}{RGB}{217, 224, 33}
\definecolor{l2yellow}{RGB}{216, 177, 64}
\definecolor{lgrey}{RGB}{179, 179, 179}
\definecolor{indigo}{rgb}{0.29, 0.0, 0.51}  
\usepackage[colorlinks, urlcolor=indigo, linkcolor=indigo, citecolor=indigo]{hyperref}

\usepackage[hcentering, vcentering, total={5.9in, 8.2in}]{geometry}  

\usepackage{tikz}
\usetikzlibrary{trees}
\usetikzlibrary{arrows}
\usepackage[all,cmtip]{xy}
 
\theoremstyle{plain}
\newtheorem{theorem}{Theorem}
\newtheorem{corollary}[theorem]{Corollary}
\newtheorem{proposition}[theorem]{Proposition}

\newtheorem{question}[theorem]{Question}
\newtheorem{conjecture}[theorem]{Conjecture}

\theoremstyle{definition}
\newtheorem{definition}[theorem]{Definition}

\theoremstyle{remark}
\newtheorem{remark}[theorem]{Remark}

\numberwithin{theorem}{section}

\usepackage{thmtools}
\usepackage{thm-restate}

\newcommand{\dfn}[1]{{\em #1}}        

\newcommand{\A}{\mathcal{A}}           
\newcommand{\FF}{\mathcal{F}}           
\newcommand{\W}{\mathcal{W}}           

\DeclareMathOperator{\bd}{\partial}   

\makeatletter
\newcommand*\bigcdot{\mathpalette\bigcdot@{0.6}}
\newcommand*\bigcdot@[2]{\mathbin{\vcenter{\hbox{\scalebox{#2}{$\m@th#1\bullet$}}}}}
\makeatother

\DeclareMathOperator{\HFhat}{\widehat{\mathit{HF}}}       
\DeclareMathOperator{\CFhat}{\widehat{\mathit{CF}}}       
\DeclareMathOperator{\HFK}{{\mathit{HFK}}}               
\DeclareMathOperator{\CFA}{\widehat{\mathit{CFA}}}        
\DeclareMathOperator{\CFD}{\widehat{\mathit{CFD}}}        

\DeclareMathOperator{\SFC}{\mathit{SFC}}       
\DeclareMathOperator{\SFH}{\mathit{SFH}}       
\newcommand{\EH}{\mathit{EH}}       

\DeclareMathOperator{\BSA}{\widehat{\mathit{BSA}}}        
\DeclareMathOperator{\BSD}{\widehat{\mathit{BSD}}}        
\DeclareMathOperator{\BSDD}{\widehat{\mathit{BSDD}}}        

\DeclareMathOperator{\AZ}{\mathsf{AZ}}          

\begin{document}

\title{Bordered contact invariants and half Giroux torsion}

\author{Hyunki Min}
\author{Konstantinos Varvarezos}

\address{Department of Mathematics \\ University of California \\ 520 Portola Plaza\\ Los Angeles, CA, United States, 90095}
\email{hkmin27@math.ucla.edu}
\email{kkv@math.ucla.edu}


\begin{abstract}
  We show that there exist infinitely many closed contact 3-manifolds containing half Giroux torsion along a separating torus whose contact invariants do not vanish. This provides counterexamples to Ghiggini's conjecture and suggests that separating half Giroux torsion may not obstruct symplectic fillability. The main tools are the bordered contact invariants recently developed by the authors and the innermost contact structures on knot complements. We also show that there exists a closed contact 3-manifold containing convex torsion along a separating torus with non-vanishing contact invariant, which implies that $2\pi$ is the minimal amount of twisting necessary to ensure vanishing of the contact invariant.     
\end{abstract}


\maketitle

\section{Introduction}\label{sec:intro}

Giroux torsion is an important class of contact structures on a neighborhood of a torus that obstructs symplectic fillability \cite{Gay:torsion,Wendl:torsion} and gives rise to infinitely many tight contact structures \cite{CGH:finite, HKM:infinite}. Furthermore, it is well known that the contact invariant in Heegaard Floer homology vanishes in the presence of Giroux torsion \cite{GHV:torsion, LS:torsion}. 

Giroux torsion along a separating torus even obstructs weak symplectic fillability \cite{GH:torsion}. Since it is significantly stronger than the non-separating case, it has been expected that even half Giroux torsion along a separating torus may obstruct symplectic fillability and could be a minimal object that does so. Along these lines, Ghiggini posed the following conjecture, which implies that separating half Giroux torsion obstructs symplectic fillability.

\begin{conjecture}[Ghiggini]\label{conjecture}
  Let $(Y, \xi)$ be a closed contact 3-manifold with half Giroux torsion along a separating torus. Then its contact invariant $\widehat{c}(\xi) \in \widehat{\mathit{HF}}(-Y)$ vanishes.
\end{conjecture}

For separating half convex torsion, which is a slight modification of half Giroux torsion (see Definition~\ref{def:convex-torsion}), Shah and Simone \cite{Shah:Seifert,SS:plumbed} showed that it does not obstruct symplectic fillability. They constructed Stein fillable contact structures on some Seifert fibered spaces with four singular fibers that contain half convex torsion along a separating torus. 

In the case of half Giroux torsion along a non-separating torus, it is well known that it does not obstruct symplectic fillability. In particular, there are contact structures on torus bundles containing half Giroux torsion for which the contact invariant does not vanish, and which admit symplectic fillings; see \cite{DG:fillability,DL:fillability, Eliashberg:T3,MRW:spinal} for example.

In the case of manifolds with boundary, Eliashberg's result on Stein fillable $T^3$ \cite{Eliashberg:T3} yields a counterexample to the conjecture: there exists a $T^2 \times I$ layer containing half Giroux torsion with non-vanishing contact invariant in sutured Floer homology. More recently, Cavallo and Matkovi\v{c} \cite{CM:torsion} gave further examples---contact structures on positive $(2,2n+1)$-torus knot complements that contain half Giroux torsion along a boundary-parallel torus. 

The conjecture has been difficult to approach in the case of closed manifolds. One reason is that it is challenging to find potential counterexamples. For instance, if the torus is not essential, the existence of half Giroux torsion in a closed manifold implies that the contact structure is overtwisted, and hence the contact invariant vanishes. Another reason is that computing the contact invariant traditionally relies on a contact surgery diagram or an open book decomposition. However, such structures become quite complicated in the presence of half Giroux torsion and make the computation of the contact invariant more difficult.

In this paper, we circumvents these difficulties by utilizing two new developments. First, we use special contact structures on knot complements, called \dfn{innermost contact structures}, to construct potential counterexamples; these are minimal contact structures on knot complements, and their precise definition can be found in Section~\ref{subsec:innermost}. Second, we utilize the bordered contact invariants recently developed by the authors \cite{MV:bordered}, which enable us to apply computational tools from bordered Floer homology to compute the contact invariants of the examples we construct. As a result, we obtain the following theorem, which provides infinitely many counterexamples to Conjecture~\ref{conjecture}.

\begin{theorem}\label{thm:half}
  For any pair of relatively prime integers $(p,q)$ with $pq < 0$, there exists a closed contact 3-manifold $(Y_{p,q},\xi_{p,q})$ with half Giroux torsion along a separating torus such that $\widehat{c}(\xi_{p,q}) \neq 0$.
\end{theorem}

Here, $Y_{p,q}$ denotes the \dfn{splice} of the left-handed trefoil and the $(p,q)$-torus knot, i.e., the manifold obtained by gluing their complements by identifying the meridian of the trefoil with the Seifert longitude of the torus knot and vice versa.

\begin{remark}\label{rmk}
  In the forthcoming paper \cite{MV:rank}, we show that the complement of every negative L-space knot admits an innermost contact structure. Moreover, the bordered invariants of L-space knot complements are structurally similar to those of torus knot complements \cite{LOT:bordered, OSz:lens, HFCC}. Thus, the same argument shows that the splice of the left-handed trefoil and any negative L-space knot admits a contact structure containing separating half Giroux torsion with non-vanishing contact invariant. 
\end{remark}

We can also ask for the minimal amount of twisting along a separating torus in a closed contact manifold that ensures the contact invariant vanishes. Recall that a Giroux torsion layer corresponds to $2\pi$-twisting, while a half Giroux torsion layer corresponds to $\pi$-twisting; see Section~\ref{subsec:torsion} for more details. To this end, we consider a contact manifold containing \dfn{convex torsion} (see Definition~\ref{def:convex-torsion}). 

\begin{theorem}\label{thm:convex}
  There exists a closed contact manifold $(Y,\xi)$ with convex torsion along a separating torus such that $\widehat{c}(\xi) \neq 0$.
\end{theorem}

Here, $Y$ denotes the double of the figure-eight knot complement. Combined with Proposition~\ref{prop:embed-torsion} stating that convex torsion contains a $(2\pi - \epsilon)$-twisting $T^2\times I$ layer with pre-Lagrangian boundaries for any small $\epsilon > 0$, we obtain the following corollary.

\begin{corollary}\label{cor:minimal}
  The minimal amount of twisting along a separating torus necessary to ensure that the contact invariant vanishes is $2\pi$.
\end{corollary}

Theorem~\ref{thm:convex} was already known for non-separating tori. For instance, in the Stein fillable contact structure on $T^3$, a Giroux torsion layer does not embed, while a convex torsion layer does \cite{Honda:classification2}.  

\begin{remark}\label{rmk:compare}
  Sometimes the terms \dfn{convex torsion} and \dfn{Giroux torsion} are used interchangeably. However, Theorem~\ref{thm:convex} shows that they are not equivalent. The difference arises from the fact that a convex torus cannot always be perturbed into a pre-Lagrangian torus. See Section \ref{subsec:torsion} for further discussion of the differences.
\end{remark}

We close the introduction with some questions for future work. First, by Theorems~\ref{thm:half} and \ref{thm:convex}, it is now conceivable that half Giroux torsion along a separating torus may not obstruct symplectic fillability. This leads to the following question. 

\begin{question}
  Are $(Y_{p,q},\xi_{p,q})$ and $(Y,\xi)$ symplectically fillable? 
\end{question}

To investigate this question, we may consider algebraic torsion \cite{LW:torsion} or filtered contact invariant \cite{KMVW:filtration}.

\begin{question}
  Compute algebraic torsion or filtered contact invariant of $(Y_{p,q},\xi_{p,q})$ and $(Y,\xi)$.
\end{question}

According to Remark~\ref{rmk}, we can construct a counterexample to Conjecture~\ref{conjecture} on the splice of the left-handed trefoil and any negative L-space knot. In general, we expect the computation to be analogous in the case where the left-handed trefoil is replaced with any negative L-space knot.

\begin{question}
  Does the splice of every pair of negative L-space knots admit a contact structure containing separating half Giroux torsion with non-vanishing contact invariant?
\end{question}

One may also ask if the construction in Theorem~\ref{thm:convex} can be generalized to an infinite family of examples. In particular, complements of even twist knots have similar bordered invariants to the complement of the figure-eight knot and according to upcoming work in \cite{MV:rank}, also admit innermost contact structures.

\begin{question}
  Does the double of every even twist knot complement admit a contact structure containing separating convex torsion with non-vanishing contact invariant?
\end{question}

By \cite[Propostion 2.5]{Martelli:Seifert}, none of the manifolds in Theorems~\ref{thm:half} and \ref{thm:convex} are Seifert fibered spaces. After several unsuccessful attempts to find counterexamples on Seifert fibered spaces, we pose the following question: 

\begin{question}
  Let $(M, \xi)$ be a closed contact Seifert fibered space with convex torsion (or half Giroux torsion) along a separating torus. Does its contact invariant $\widehat{c}(\xi)$ always vanish?
\end{question}

\section*{Acknowledgement}

The authors thank Joan Licata for introducing them to the conjecture.

\section{Background}\label{sec:background}

In this section, we recall the relevant background material we will need from contact geometry and bordered Floer homology. For more details, the reader is referred to \cite{Etnyre:convex, Geiges:book} for contact geometry and \cite{HRW:immersed, LOT:bordered} for bordered Floer homology.

\subsection{Giroux and convex torsion}\label{subsec:torsion}
We first recall the definition of Giroux torsion.
\begin{definition}
  For $n \in \frac12\mathbb{N}$, the following contact manifold
  \[
    (T^2 \times [0,1], \ker(\cos(2\pi nz)dx + \sin(2\pi nz)dy))
  \]
  is called a \dfn{Giroux $n$-torsion layer} or a \dfn{$2n\pi$-twisting layer}. We say a contact manifold $(Y,\xi)$ \dfn{contains Giroux torsion} (resp. \emph{half Giroux torsion}) if it contains an embedded Giroux $n$-torsion layer for some $n \geq 1$ (resp. $n \geq \frac12$). 
\end{definition}

Once we fix a basis of $H_1(T^2)$, we can extend the previous definition to any rational amount of twisting. 

\begin{definition} \label{def:twisting-layer}
  For $r \in \mathbb{Q}^+$, the following contact manifold
  \[
    (T^2 \times [0,1], \ker(\cos(2\pi rz)dx + \sin(2\pi rz)dy))
  \] 
  is called \dfn{$2r\pi$-twisting layer}. 
\end{definition}

The characteristic foliations of the boundary tori of a Giroux torsion layer are pre-Lagrangian. We can also consider an analogous version in which the boundaries are convex, which is called \dfn{convex torsion}.

\begin{definition} \label{def:convex-torsion}
  For $n \in \frac12\mathbb{N}$ and a small $\epsilon >0$, consider the following contact manifold: 
  \[
    (T^2 \times [-\epsilon,1+\epsilon], \ker(\cos(2\pi nz)dx + \sin(2\pi nz)dy)).
  \]
  Perturb $T^2 \times \{0\}$ and $T^2 \times \{1\}$, fixing a single leaf of their characteristic foliations, to make them convex; denote the resulting convex tori by $T_0$ and $T_1$, respectively. The contact manifold bounded by $T_0$ and $T_1$ is called a \dfn{convex $n$-torsion layer}. We say a contact manifold $(Y,\xi)$ \dfn{contains convex torsion} (resp. \emph{half convex torsion}) if it contains an embedded convex $n$-torsion layer for some $n \geq 1$ (resp. $n \geq \tfrac12$).
\end{definition}

\begin{remark}
  The $\epsilon$-extension of $T^2 \times [0,1]$ is necessary in the above definition. Since we fix a leaf of the characteristic foliation, the surface can be viewed as having boundary. To make it convex near the boundary, we use a standard annulus model which requires a two-sided collar neighborhood; see \cite[Theorem 2.23]{Etnyre:convex}.

\end{remark}

The following proposition is a direct consequence of the previous definitions.

\begin{proposition}\label{prop:embed-torsion}
  For any $n, m \in \frac12\mathbb{N}$ with $n < m$, a Giroux $n$-torsion layer  embeds in a convex $m$-torsion layer. Also, a convex $n$-torsion layer embeds in a Giroux $m$-torsion layer.

  Moreover, for any $n\in \frac12\mathbb{N}$ and $r \in \mathbb{Q}$ such that $0 < r < n$, a $2r\pi$-twisting layer embeds in a convex $n$-torsion layer.
\end{proposition}

As mentioned in the introduction, the contact invariant of a contact manifold containing Giroux torsion vanishes.

\begin{theorem}[Lisca--Stipsicz \cite{LS:torsion}, Ghiggini--Honda--Van-Horn-Morris \cite{GHV:torsion}]
  Let $(Y,\xi)$ be a closed contact 3-manifold containing Giroux torsion. Then the contact invariant $\widehat{c}(\xi) \in \widehat{\mathit{HF}}(-Y)$ vanishes.
\end{theorem}

There are subtle differences between Giroux torsion and convex torsion. In particular, Proposition~\ref{prop:embed-torsion} does not hold when $n=m$. In particular, there are contact $3$-manifolds that contain convex torsion but not Giroux torsion. For example, the contact structure on $T^3$ 
\[
  \xi = \ker(\cos(2\pi z)dx + \sin(2\pi z)dy) 
\]
is Stein fillable, and hence its contact invariant is non-vanishing, which implies that it does not contain Giroux torsion. However, a convex $1$-torsion layer embeds in $(T^3, \xi)$. The embedding can be seen as follows: take a torus $T^2 \times \{0\}$ in $T^3$ and make it convex. Then take an $I$-invariant neighborhood of this torus and remove it from $T^3$. The remaining part is a convex $1$-torsion layer. This implies that convex $1$-torsion does not obstruct symplectic fillability, while Giroux $1$-torsion does \cite{Gay:torsion,Wendl:torsion}. The following result is also immediate.

\begin{proposition}
  Let $(T^2\times[0,1],\xi_1)$ be a convex $1$-torsion layer. Then its sutured contact invariant $\EH(\xi_1)$ is non-vanishing. 
\end{proposition}

\subsection{Bordered contact invariants}\label{subsec:bordered}
In this section, we briefly review some relevant facts about bordered contact invariants from \cite{MV:bordered}. Let $(Y,\FF,\Gamma)$ be a bordered sutured 3-manifold, where $\FF\subset \bd Y$ is parametrized, and $\Gamma$ is a set of dividing curves on $\bd Y \setminus \FF$. Also let $\A=\A(\FF)$ be the strand algebra associated to the surface $\FF.$  To $(Y,\FF,\Gamma)$ we may associate a type-A module $\BSA(Y)$, which is an $\A_{\infty}$ module over $\A$, as well as a type-D structure $\BSD(Y)$. There is an algebraic pairing of bordered modules that recovers the sutured Floer homology of a pair of bordered sutured manifolds glued along their common parametrized boundaries: if $(Y_1,\FF, \Gamma_1)$ and $(Y_2,-\FF, \Gamma_2)$ are bordered sutured 3-manifolds, then
\[
\SFC(Y_1 \cup_{\FF} Y_2,\Gamma_1\cup\Gamma_2) \cong \BSA(Y_1)\boxtimes \BSD(Y_2)\cong \BSD(Y_2)\boxtimes\BSA(Y_2).
\]

For a bordered sutured 3-manifold $(Y,\FF,\Gamma)$ and any idempotent $\iota \in \A$, there exists a bordered sutured manifold $\W_\iota=(\FF\times [0,1], -\FF,\Gamma_{\iota})$ called the  \emph{elementary cap} associated to $\iota$. We have the following natural identification:
\begin{equation}\label{eq:cap}
\BSA(Y)\cdot \iota = \BSA(Y)\boxtimes \BSD(\W_{\iota}) \cong \SFC(Y,\Gamma\cup\Gamma_{\iota})
\end{equation}

Now suppose $\xi$ is a contact structure on $Y$ such that for some idempotent $\iota \in \A,$ the contact structure induces dividing curve $\Gamma \cup \Gamma_{\iota}$. Then there is a well-defined contact invariant $c_A(\xi)\in\BSA(Y)$, characterized by 
\[
 c_A(\xi)=c_A(\xi)\cdot \iota,\quad [c_A(\xi)\cdot \iota] = EH(\xi)\in \SFH(-Y,-\Gamma \cup -\Gamma_{\iota}) 
\]
under the identification \eqref{eq:cap}.

Given $c_A(\xi)$ we may also define $c_D(\xi)=c_A(\xi)\boxtimes \iota^{\vee} \in \BSA(Y)\boxtimes \BSDD(\AZ).$ Here, $\AZ$ is a specific Heegaard Diagram depending only on $\FF,$ and $\iota^{\vee}$ is a specific generator corresponding to the idempotent $\iota.$ In the case where $\FF$ is a once-punctured torus, an explicit model for $\BSDD(\AZ)$ was presented in \cite[Section 4.1]{MV:bordered}: it is generated by $\iota_0^{\vee}, \iota_1^{\vee}, \rho_{1}^{\vee}, \rho_{2}^{\vee}, \rho_{3}^{\vee}, \rho_{12}^{\vee}, \rho_{23}^{\vee}, \rho_{123}^{\vee}$ and the differential is given by 
\begin{align*}
  \delta^1\left(\rho_{123}^{\vee}\right) &= \rho_3 \otimes \rho_{12}^{\vee} + \rho_{23}^{\vee} \otimes \rho_1, & \delta^1\left(\rho_{23}^{\vee}\right) &= \rho_3 \otimes \rho_{2}^{\vee} + \rho_{3}^{\vee} \otimes \rho_2, \\
  \delta^1\left(\rho_{12}^{\vee}\right) &= \rho_2 \otimes \rho_{1}^{\vee} + \rho_{2}^{\vee} \otimes \rho_1, & \delta^1\left(\rho_{3}^{\vee}\right) &= \rho_3 \otimes \iota_0^{\vee} + \iota_1^{\vee} \otimes \rho_3, \\
  \delta^1\left(\rho_{1}^{\vee}\right) &= \rho_1 \otimes \iota_0^{\vee} + \iota_1^{\vee} \otimes \rho_1, & \delta^1\left(\rho_{2}^{\vee}\right) &= \rho_2 \otimes \iota_1^{\vee} + \iota_0^{\vee} \otimes \rho_2,
\end{align*}
where $\{\iota_0, \iota_1, \rho_1, \rho_2, \rho_3, \rho_{12}, \rho_{23}, \rho_{123}\}$ forms the torus algebra. 

We have the following pairing theorem, which allows us to recover the sutured contact invariant by pairing type-A and type-D invariants. 

\begin{theorem}[Theorem 1.2 of \cite{MV:bordered}]\label{thm:contactPairing}
If $(Y_1,\FF, \Gamma_1)$ and $(Y_2,-\FF, \Gamma_2)$ are bordered sutured 3-manifolds with contact structures $\xi_1$ and $\xi_2,$ respectively (compatible with the parametrization of $\FF$), then
\[
[c_A(\xi_1)\boxtimes c_D(\xi_2)] = \EH(\xi_1 \cup \xi_2) \in \SFH(Y_1 \cup_{\FF} Y_2,\Gamma_1\cup\Gamma_2)
\]
\end{theorem}

For a bordered sutured 3-manifold with torus boundary $(Y,\FF, \Gamma),$ if $\FF=T^2 \setminus D^2$ and $\Gamma$ is a trivial arc in $D^2,$ unlinked with respect to the parametrizing arcs, we can identify $(Y,\FF,\Gamma)$ with a bordered manifold $(Y,\FF)$ (treating the disk as a point). Thus the previous statements still hold after replacing $\BSA$ with $\CFA$, $\BSD$ with $\CFD$, $\SFC$ with $\CFhat$ and $EH$ with $\widehat{c}$. In this paper, we only consider closed manifolds and manifolds with torus boundary, so we assume this condition on $(Y,\FF,\Gamma)$ whenever a torus boundary is involved. 

For a given parametrized surface $\FF,$ Mathews \cite{Mathews:contactCategory} has demonstrated a correspondence between homology classes in $\A(\FF)$ and tight contact structures on $\FF \times [0,1]$. For a cycle $a\in\A(\FF)$, let us denote by $\xi_a$ the corresponding contact structure on $\FF \times [0,1]$.

\begin{theorem}[Corollary~1.9 of \cite{MV:bordered}]\label{thm:Ainf_bypass}
Let $(Y,\FF,\Gamma)$ be a bordered sutured manifold with a compatible contact structure $\xi$. Then for every cycle $a\in\A(\FF)$, we have
\begin{align*}
c_A(\xi\cup \xi_a) = m_2(c_A(\xi),a).
\end{align*}
\end{theorem}

In the case of torus boundary, we have explicit description of the correspondence between cycles in $\A(\FF)$ and tight contact structures on $\FF \times [0,1]$. 

\begin{theorem}[Theorem 1.7 of \cite{MV:bordered}]\label{thm:basicslice} 
  In the case of torus boundary, each idempotent corresponds to a dividing set on the boundary torus $\Gamma_0$ and $\Gamma_1$. Each $\Gamma_i$ consists of two dividing curves parallel to one parametrizing arc of $\FF$. Generators of $\A(\FF)$ correspond to tight contact structures on $T^2 \times [0,1]$ as follows:

  \begin{itemize}
    \item $\rho_1$, $\rho_3$: basic slices with boundary dividing sets $\Gamma_0$ and $\Gamma_1$ (they have different signs).
    \item $\rho_2$: a basic slice with boundary dividing sets $\Gamma_1$ and $\Gamma_0$.
    \item $\rho_{12}$: a half convex torsion layer with boundary dividing set $\Gamma_0$.
    \item $\rho_{23}$: a half convex torsion layer with boundary dividing set $\Gamma_1$.
    \item $\rho_{123}$: a union of a basic slice and a half convex torsion layer ($3\pi/2$-layer). The boundary dividing sets are $\Gamma_0$ and $\Gamma_1$. 
  \end{itemize}
\end{theorem}

The following proposition is an immediate consequence of Theorem~\ref{thm:basicslice} and the definitions in Section~\ref{subsec:torsion}.

\begin{proposition}\label{prop:123torsion}
  The contact structure corresponding to $\rho_{123}$ contains half Giroux torsion.
\end{proposition}

\begin{figure}
\centering
\begin{subfigure}{.5\textwidth}
  \centering
  \footnotesize
  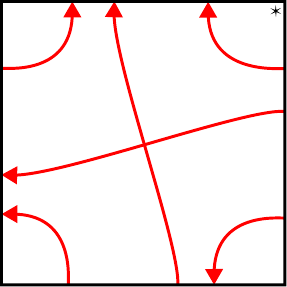
  \caption{Immersed curve type A operations}
  \end{subfigure}\hfill
  \begin{subfigure}{.5\textwidth}
  \centering
  \footnotesize
  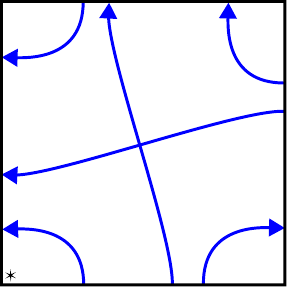
  \caption{Immersed curve type D operations}
\end{subfigure}\\
\vspace*{.5cm}
\begin{subfigure}{.5\textwidth}
  \centering
  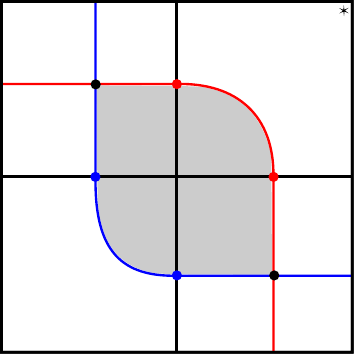
  \caption{pairing immersed curves}
  \label{fig:Imm_pairing}
\end{subfigure}
\caption{Immersed curve interpretation of bordered invariants and their pairing}
\end{figure}

In the case of torus boundary, the work of Hanselman--Rasmussen--Watson interprets the bordered invariants and the pairing theorem in terms of immersed curves in the punctured torus $T^2\setminus \{*\}$ \cite{HRW:immersed}. Indeed, the punctured torus is canonically identified with $\FF$ so that the generators of the bordered module correspond to intersection points between the immersed curves and the parametrizing arcs. The pairing of type-A and type-D modules  may be performed by homotoping the associated curves to lie in the upper right and lower left of the torus (in the usual square fundamental domain, thought of as the complement of the parametrizing arcs which are the edges of the square) except for purely vertical and horizontal arcs. The intersection points between the type-A and type-D arcs correspond to generators of the paired complex. The differential involves counting immersed bigons between generators. See Figure~\ref{fig:Imm_pairing}.

In the case of the complement of a knot $K\subset S^3$, the immersed curve invariant encodes several invariants from Knot Floer Homology. There is a standard lift of the immersed curve invariant of $S^3\setminus K$ to the infinitely punctured cylinder (see \cite[Section 6]{HRW:immersed}), such that the meridian edge lifts to the vertical unbounded line, whereas the Seifert longitude is the horizontal closed curve. 
The intersection points of the immersed curve with the vertical line correspond to generators of $\widehat{\HFK}(S^3,K)$, and the height of each intersection corresponds to the Alexander grading of its corresponding generator.

\subsection{Innermost contact structures on knot complements}\label{subsec:innermost}
We review \dfn{innermost contact structures} on knot complements. These were first introduced by Conway and the first author for the figure-eight knot \cite{CM:figure8} in the classification of tight contact structures, and were later extended to torus knots \cite{EMM:non-loose} and to certain knots in lens spaces \cite{MN:whitehead}.

Let $K$ be a knot in a closed 3-manifold $Y$ and $C = Y \setminus \nu(K)$ denote the complement of $K$. We define an innermost contact structure on $C$ as follows.

\begin{definition}
  We say a contact structure $\xi$ on $C$ is \dfn{innermost} if
  \begin{enumerate}
    \item $\partial C$ is convex with two meridional dividing curves, and
    \item $\partial C$ does not thicken, i.e., no contact structure $\xi'$ on $C$ that induces non-meridional dividing curves on $\partial C$ embeds into $\xi$.
  \end{enumerate}
\end{definition}

An innermost contact structure can be considered as a minimal element in the set of contact structures on $C$. In general, it is hard to find innermost contact structures, since they do not embed in tight contact structures on $Y$, but only into overtwisted ones. Moreover, not every knot admits an innermost contact structure. For instance, the complements of the unknot and the right-handed trefoil do not admit one. 

In this paper, we utilize innermost contact structures of the figure-eight knot and negative torus knots. In particular, these contact structures have two important properties: their contact invariants are non-vanishing in sutured Floer homology, and their Alexander grading is $g(K) - 1$ where $g(K)$ is the Seifert genus of the knot. 

\begin{proposition}[Conway--Min \cite{CM:figure8}]\label{Prop:fig8_inner}
  The figure-eight knot complement admits an innermost contact structure $\xi^{\textit{in}}_8$ satisfying the following properties:
  \begin{enumerate}
    \item Its sutured contact invariant $\EH(\xi^{\textit{in}}_8)$ is non-vanishing.
    \item The Alexander grading of $\EH(\xi^{\textit{in}}_8)$ is $0$.
    \item $\xi^{\textit{in}}_8$ embeds into the contact structure on $S^3$ supported by the figure-eight knot, which is overtwisted.
    \item After attaching a half convex torsion layer to the boundary of $\xi^{\textit{in}}_8$, the resulting contact structure still has a non-vanishing contact invariant.
  \end{enumerate} 
\end{proposition}

\begin{proposition}[Etnyre--Min--Mukherjee \cite{EMM:non-loose}]\label{Prop:Tpq_inner}
  The complement of any negative torus knot $T_{p,q}$ admits an innermost contact structure $\xi^{\textit{in}}_{p,q}$ satisfying
  \begin{enumerate}
    \item Its sutured contact invariant $\EH(\xi^{\textit{in}}_{p,q})$ is non-vanishing.
    \item The Alexander grading of $\EH(\xi^{\textit{in}}_{p,q})$ is $g(T_{p,q})-1$. 
    \item $\xi^{\textit{in}}_{p,q}$ embeds into the contact structure on $S^3$ supported by $T_{p,q}$, which is overtwisted.
  \end{enumerate} 
\end{proposition}

We briefly describe how to obtain the innermost contact structures of the left-handed trefoil and the figure-eight knot. Since both are genus-one fibered knots, $0$-surgery yields torus bundles over $S^1$. Each of these torus bundles admits a Stein fillable contact structure. Now take the surgery dual knot in each Stein fillable torus bundle and make it Legendrian with maximal twisting number. By removing a standard neighborhood of this surgery dual knot, we obtain an innermost contact structure on the knot complement. 

\section{Proof of Theorems}\label{sec:proof}

In this section, we prove Theorems~\ref{thm:half} and \ref{thm:convex}. We begin with the proof of Theorem~\ref{thm:half}, which uses complements of negative torus knots. Then we prove Theorem~\ref{thm:convex} using the figure-eight knot complement. The brief outline of the proof is as follows: take two innermost contact structures on knot complements and insert a $T^2 \times I$ layer that contains half Giroux torsion or convex torsion. We then apply the pairing theorem to show that the contact invariant of the resulting closed contact manifold is non-vanishing.

\begin{figure}
\centering
\footnotesize
	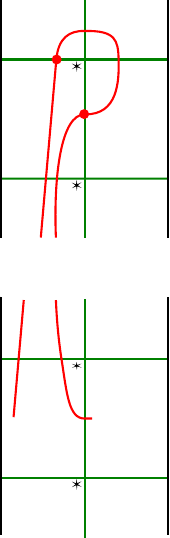
\caption{Part of the (lift to the cylinder of the)  immersed curve of the complement of a positive torus knot}
\label{fig:Tpq-immersed}
\end{figure}

\subsection{The torus knot complements}
We first identify the bordered contact invariants of innermost contact structures on torus knot complements. The bordered modules of torus knot complements are well-understood; for instance, they can be derived from the corresponding knot Floer complexes via the procedure described in \cite[Chapter 11]{LOT:bordered} (the immersed curve version of this procedure is outlined in \cite[Section 4]{HRW:immersed2}). We focus on the bordered modules of positive torus knots, since we work with contact structures on complements of negative torus knots, and the contact invariant lives in the bordered module of the manifold with reversed orientation. 

\begin{figure}
\centering
\raisebox{5.5cm}{
\begin{minipage}{.3\textwidth}
\begin{subfigure}{\textwidth}
  \footnotesize
  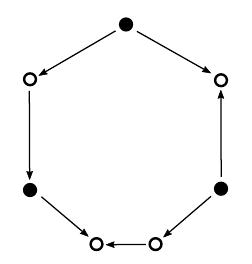
  \caption{A (left) type A module  for the Seifert-framed right-handed trefoil complement}
  \label{fig:trefoil-A}
\end{subfigure}\\
\vspace{.5 cm}
\begin{subfigure}{\textwidth}
  \addtocounter{subfigure}{1}
  \footnotesize
  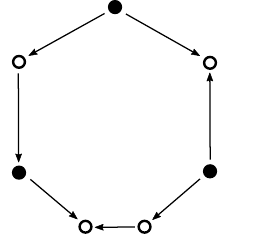
  \caption{Simplified type D structure for the Seifert-framed right-handed trefoil complement}
  \label{fig:trefoil-D}
\end{subfigure}
\end{minipage}
}
\hfill
\begin{subfigure}{.65\textwidth}
  \addtocounter{subfigure}{-2}
  \centering
  \footnotesize
  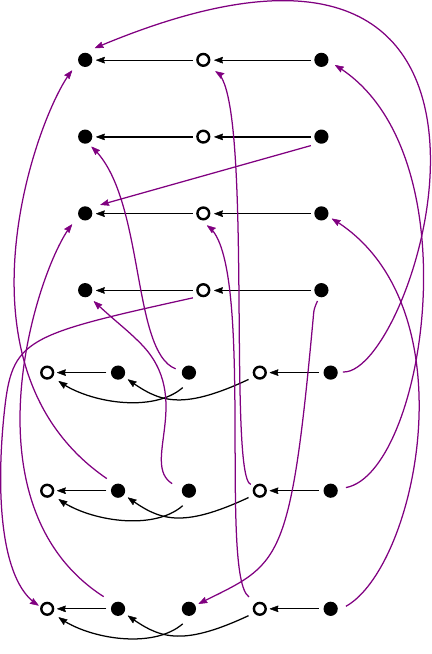
  \caption{Type D structure for trefoil complement obtained by pairing type-A module with $\BSDD(\AZ)$. Here, $M = c_A(\xi^{\textit{in}}_{-2,3})$}
  \label{fig:trefoil-DF}
\end{subfigure}
\caption{Bordered contact invariants for the trefoil complement}
\end{figure}

Let $T_{p,q}$ be a negative torus knot where $(p,q)$ is a pair of relatively prime integers with $p < 0 < q$. The lifted immersed curve invariant for a positive torus knot $T_{-p,q}$ consists of a single embedded curve whose behavior near the maximal Alexander grading is as shown in Figure~\ref{fig:Tpq-immersed}. There is a unique generator in Alexander grading $g(T_{p,q})-1$, so by Proposition~\ref{Prop:Tpq_inner}, this must be the contact invariant $c_A(\xi^{\textit{in}}_{p,q})$ for the innermost contact structure $\xi^{\textit{in}}_{p,q}$. From the structure of the immersed curve, we see that $m_2(c_A(\xi^{\textit{in}}_{p,q}),\rho_{123})$ is a nontrivial generator. By Theorem~\ref{thm:Ainf_bypass}, it corresponds to attaching a contact structure $\xi_{123}$ corresponding to $\rho_{123}$, i.e.,
\[
 m_2(c_A(\xi^{\textit{in}}_{p,q}),\rho_{123}) = c_A(\xi^{\textit{in}}_{p,q}\cup\xi_{123}).  
\]
We denote the resulting contact structure by $\xi_{p,q}^{123}$. Since we parametrize the boundary torus with meridian and Seifert longitude, one boundary component of $\xi_{123}$ has meridional dividing curves, and the other has longitudinal dividing curves by Theorem~\ref{thm:basicslice}. Since we glue the contact structures along the torus with meridional dividing curves, we know $\xi_{p,q}^{123}$ induces longitudinal dividing curves on the boundary. By Proposition~\ref{prop:123torsion}, the contact structure $\xi_{123}$ contains a half Giroux torsion layer, so $\xi_{p,q}^{123}$ contains a half Giroux torsion layer along a boundary-parallel torus. 

\begin{figure}
\footnotesize
  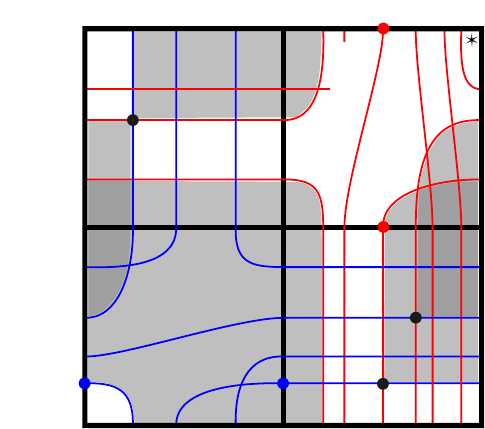
  \caption{Immersed curve computation of the contact invariant for $Y_{p,q}$}
  \label{fig:T-splice}
\end{figure}

In the case of the trefoil, Figure~\ref{fig:trefoil-A} shows a type-A module for the Seifert-framed right-handed trefoil complement. There is a unique generator corresponding to Alexander grading $0=g(T_{2,3})-1$, so Proposition~\ref{Prop:Tpq_inner} again implies that it should be $\xi^{\textit{in}}_{-2,3}$. To find the corresponding type-D contact invariant, we pair the type-A module with the DD-bimodule $\BSDD(\AZ)$. Figure~\ref{fig:trefoil-DF} shows the result of this pairing. The purple arrows are all pure differentials and may be cancelled. The resulting simplified type-D module is shown in Figure~\ref{fig:trefoil-D}, and the type-D invariant is $\iota^{\vee}_0 c_A(\xi^{\textit{in}}_{-2,3}).$

\begin{remark}
  Instead of using the splice, if we identify the meridian (resp. longitude) of the trefoil with the meridian (resp. longitude) of the torus knot, then we cannot insert $\rho_{123}$, as the dividng slopes do not match. Instead, we can insert a half convex torsion layer corresponding to $\rho_{12}$ and the resulting contact structure may not contain half Giroux torsion. If we insert a full convex torsion layer, then the resulting contact structure will have vanishing contact invariant.
\end{remark}

\subsection{Proof of Theorem~\ref{thm:half}}
Given a pair of relatively prime integers $(p,q)$ with $pq<0$, we define $Y_{p,q}$ to be the splice of the left-handed trefoil and the torus knot $T_{p,q}$; that is,
\[
Y_{p,q} = S^3\setminus \nu(T_{p,q}) \cup_{\phi} S^3\setminus \nu(T_{-2,3})
\]
where $\phi$ identifies a meridian of the trefoil with a Seifert longitude of $T_{p,q}$, and vice versa. The contact structure $\xi_{p,q}$ on $Y_{p,q}$ is obtained by gluing the contact structures $\xi^{123}_{p,q}$ and $\xi^{\textit{in}}_{-2,3}$. This gluing is well-defined since  $\xi^{\textit{in}}_{-2,3}$ induces meridional dividing curves on the boundary, while $\xi^{123}_{p,q}$ induces longitudinal dividing curves, so $\phi$ maps one dividing set to the other. Since $(p,q)$ is fixed, we denote $\xi_{p,q}$ simply by $\xi$ for brevity. As discussed in the previous subsection, $\xi^{123}_{p,q}$ contains a half Giroux torsion layer along a boundary-parallel torus of $S^3 \setminus \nu(T_{p,q})$, so $\xi$ contains half Giroux torsion along a separating torus of $Y_{p,q}$.

According to Theorem~\ref{thm:contactPairing}, $c(\xi) := c_A(\xi^{123}_{p,q}) \boxtimes c_D(\xi^{\textit{in}}_{-2,3})$ represents the contact class of $\xi$, that is, $[c(\xi)] = \widehat{c}(\xi) \in \HFhat(-Y_{p,q})$. We will show that $c(\xi)$ is nontrivial in homology. To do so, we perform the bordered paring computation using immersed curves, shown in Figure~\ref{fig:T-splice}. In the top right, we have part of the type-A immersed curve of the torus knot complement (this is simply the projection of Figure~\ref{fig:Tpq-immersed} to the torus) with the generator corresponding to $c_A(\xi_{p,q}^{123})$ marked. In the bottom left, we have the type-D module of the trefoil complement with the generator corresponding to $c_D(\xi_{-2,3}^{\textit{in}})$ marked. The resulting contact class $c(\xi)$ is  the intersection of the horizontal and vertical arcs corresponding to the two bordered contact generators. We see that there are generators $x$ and $y$, as well as an immersed bigon from $x$ to the contact class $c(\xi)$ and another from $x$ to $y$. As these are the only bigons involving these three generators, these generators generate a summand in $\CFhat(-Y_{p,q})$ with differential given by $\partial x = y+c(\xi).$ Hence, $c(\xi)$ survives in homology.

\begin{figure}
\centering
\raisebox{7cm}{
\begin{minipage}{.3\textwidth}
\begin{subfigure}{\textwidth}
  \footnotesize
  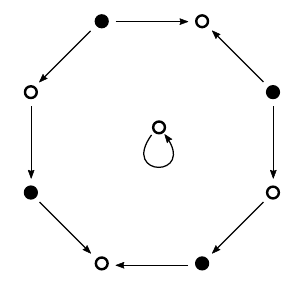
  \caption{A (left) type A module  for the Seifert-framed Figure-8 knot complement}
  \label{fig:fig8-A}
\end{subfigure}\\
\vspace{.5cm}
\begin{subfigure}{\textwidth}
  \addtocounter{subfigure}{1}
  \footnotesize
  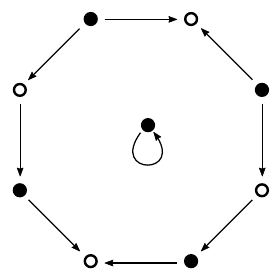
  \caption{Simplified type D structure for the Seifert-framed Figure-8 complement}
  \label{fig:fig8-D}
\end{subfigure}
\end{minipage}
}
\hfill
\begin{subfigure}{.65\textwidth}
  \addtocounter{subfigure}{-2}
  \centering
  \footnotesize
  \import{figures}{Fig8_D_full.pdf_tex}
  \caption{Type D structure for the figure-eight knot complement obtained by pairing type-A module with $\BSDD(\AZ)$}
  \label{fig:fig8-DF}
\end{subfigure}
\caption{Type A and type D modules for the figure-8 knot complement}
\label{fig:Fig8-bordered}
\end{figure}

\subsection{The figure-eight knot complement}
We first identify the bordered contact invariants of an innermost contact structure on the figure-eight knot complement. As in the torus knot case, the bordered modules for the figure-eight knot complement can be derived from the knot Floer complex. See Figures~\ref{fig:fig8-A} and \ref{fig:fig8-D}.

In the type-A module in Figure~\ref{fig:fig8-A}, the only generators of Alexander grading $0=g(4_1)-1$ are $E$ and $I$. Hence by Proposition~\ref{Prop:fig8_inner}, $c_A(\xi^{\textit{in}}_8)$ is a linear combination of these two generators. Furthermore, since the contact invariant of $\xi^{\textit{in}}_8$ is non-vanishing after attaching a half convex torsion layer by Proposition~\ref{Prop:fig8_inner}, we know $m_2(c_A(\xi^{\textit{in}}_8),\rho_1)\neq 0$. Thus $c_A(\xi^{\textit{in}}_8)$ must be one of $E$ or $E+I$.

From the type-A module in Figure~\ref{fig:fig8-A}, we can see that $m_2(E,\rho_{23})=m_2(E+I,\rho_{23})=B.$ Hence, by Theorem~\ref{thm:Ainf_bypass} and \ref{thm:basicslice}, the generator $B$ is the contact invariant for the contact structure obtained from $\xi^{\textit{in}}_8$ by attaching a half convex torsion layer corresponding to $\rho_{23}$. Thus the resulting contact structure contains half convex torsion and induces meridional dividing curves on the boundary. Since $B$ and $T$ are the unique generators in Alexander gradings $1$ and $-1$, respectively, they represent the contact invariants of conjugate contact structures (recall that conjugate contact structures can be obtained by reversing the coorienation). Hence, $T$ is also the contact invariant of a contact structure which contains a half convex torsion layer with meridional dividing curves on the boundary, and we further see that $m_2(T,\rho_{23})=I$. Let $\xi_8^{\textit{tor}}$ be the result of gluing the contact structures corresponding to $T$ and $\rho_{23}$. Then $c_A(\xi_8^{\textit{tor}}) = I$. Also, $\rho_{23}$ corresponds to a half convex torsion layer and $T$ already contains half convex torsion, so $\xi_8^{\textit{tor}}$ contains a full convex torsion layer along a boundary-parallel torus and induces meridional dividing curves on the boundary. 

To find the type-D contact invariant of $\xi_8^{\textit{in}}$, we pair the type-A module with the DD-bimodule $\BSDD(\AZ)$. Figure~\ref{fig:fig8-DF} shows the result of this pairing. The purple arrows are all pure differentials, and may be cancelled. The resulting simplified type-D module is shown in Figure~\ref{fig:fig8-D}, and the type-D invariant is $\iota^{\vee}_0 c_A(\xi^{\textit{in}}_{8})$, which is either $\iota^{\vee}_0 \boxtimes E$ or $\iota^{\vee}_0 \boxtimes (E+I)$.

\begin{remark}
  The figure-eight knot complement admits two conjugate innermost contact structures. It is known that their contact invariants are distinct. Thus both $E$ and $E+I$ correspond to the contact invariants of these innermost contact structures.
\end{remark}

\begin{figure}
\centering
\scriptsize
\begin{subfigure}{.48\textwidth}
  \centering
  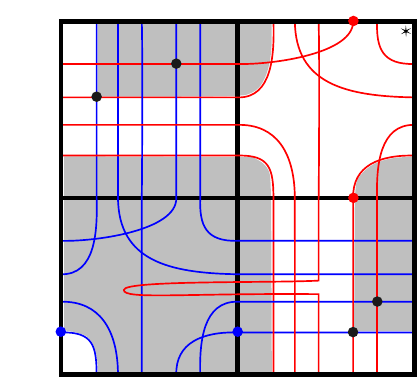
  \caption{An immersed bigon from $z$ to $c(\xi)$}
  \label{fig:Fig8-immersedA}
\end{subfigure}
\hfill
\begin{subfigure}{.48\textwidth}
  \centering
  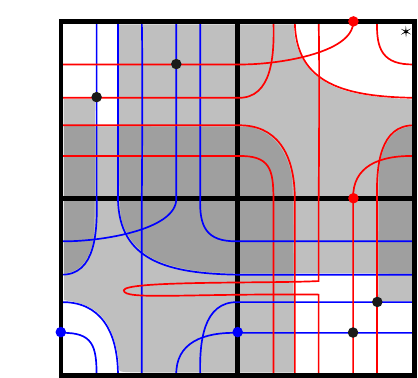
  \caption{An immersed bigon from $z$ to $y$}
  \label{fig:Fig8-immersedB}
  \end{subfigure}\\
  \vspace*{.5cm}
\begin{subfigure}{.48\textwidth}
  \centering
  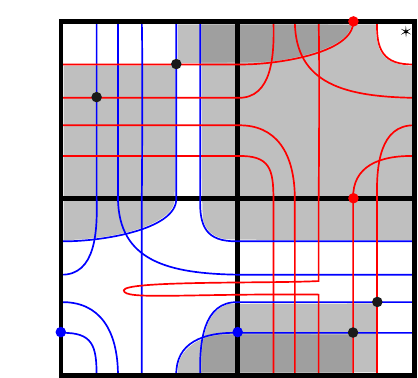
  \caption{Immersed bigons from $x$ to $y$ and from $x$ to $c(\xi)$}
  \label{fig:Fig8-immersedC}
\end{subfigure}
\caption{Immersed curve computation of contact invariant for the double of the Figure-eight knot complement}
\label{fig:Fig8-immersed}
\end{figure}

\subsection{Proof of Theorem~\ref{thm:convex}}
Consider the manifold $Y$ obtained by doubling the figure-eight knot complement:
\[
Y=S^3\setminus \nu(4_1) \cup_{\phi} S^3\setminus \nu(4_1),
\]
where $\phi$ identifies the two meridians and the two Seifert longitudes. The contact structure $\xi$ on $Y$ is obtained by gluing the contact structures $\xi^{\textit{in}}_8$ and $\xi_8^{\textit{tor}}$. This gluing is well-defined since both $\xi^{\textit{in}}_8$ and $\xi_8^{\textit{tor}}$ induce meridional dividing curves on the boundary so $\phi$ sends one dividing set to the other. As discussed in the previous subsection, $\xi_8^{\textit{tor}}$ contains a convex torsion layer along a boundary parallel-torus of $S^3\setminus \nu(4_1)$, so $\xi$ contains convex torsion along a separating torus of $Y$. 

According to Theorem~\ref{thm:contactPairing}, $c(\xi) := c_A(\xi^{tor}_{8}) \boxtimes c_D(\xi^{\textit{in}}_{8})$ represents the contact class of $\xi$, i.e., $[c(\xi)] = \widehat{c}(\xi) \in \HFhat(-Y)$. We will show that $c(\xi)$ is nontrivial in homology. To do so, we first note that since $\xi_8^{\textit{tor}}$ contains a convex torsion layer along a boundary-parallel torus, so $\xi_8^{\textit{tor}} \cup \xi_8^{\textit{tor}}$ contains a convex $2$-torsion layer and hence contains Giroux torsion by Proposition~\ref{prop:embed-torsion}. Also recall that $c_A(\xi_8^{\textit{tor}})=I$. Thus $I\boxtimes \iota_0^{\vee}\boxtimes I$ is nullhomologous as it corresponds to the contact class of $\xi_8^{\textit{tor}} \cup \xi_8^{\textit{tor}}$. Hence, $I\boxtimes \iota_0^{\vee}\boxtimes E$ and $I\boxtimes \iota_0^{\vee}\boxtimes (E+I)$ represent the same homology class, namely, 
\[
  [c_A(\xi^{tor}_{8})\boxtimes \iota_0^{\vee}\boxtimes c_A(\xi^{\textit{in}}_{8})] = [c_A(\xi^{tor}_{8}) \boxtimes c_D(\xi^{\textit{in}}_{8})] = \widehat{c}(\xi).
\] 
So it suffices to perform our calculations assuming $c_D(\xi_{in})=\iota_0^{\vee}\boxtimes E$. We show those calculations performed using immersed curves in Figure~\ref{fig:Fig8-immersed}.

In the top right of the squares in Figure~\ref{fig:Fig8-immersed}, we have part of the type-A immersed curve of the figure-eight knot complement with the generator corresponding to $c_A(\xi_8^{\textit{tor}})$ marked. Since we are gluing meridian to meridian, we reparametrize the type-A immersed curve by a $\pi/2$-rotation.  In the bottom left, we have the type-D module of the figure-eight knot complement with the generator corresponding to $c_D(\xi_{8}^{\textit{in}})$ marked.
The resulting contact class $c(\xi)$ is the intersection of the horizontal and vertical arcs corresponding to the two bordered contact generators. For admissibility reasons, the vertical type-A curve has been homotoped to intersect the vertical type-D curve; this however does not affect our computation, as the relevant generators do not belong to theses curves.

In Figure~\ref{fig:Fig8-immersed}, we see that there are generators $x$ and $z$ with bigons from each of these to the contact class $c(\xi).$ Furthermore, there are bigons from each of $x$ and $z$ to another generator, labelled $y$ in Figure~\ref{fig:Fig8-immersed}. As these are the only bigons involving these four generators, they generate a summand in $\CFhat(-Y)$ with differential given by $\partial x = y+c(\xi)=\partial z.$ Hence, $c(\xi)$ survives in homology.

\bibliography{references}
\bibliographystyle{plain}
\end{document}

%% file: figures/TypeA_labels_s.pdf_tex
\begingroup%
  \makeatletter%
  \providecommand\color[2][]{%
    \errmessage{(Inkscape) Color is used for the text in Inkscape, but the package 'color.sty' is not loaded}%
    \renewcommand\color[2][]{}%
  }%
  \providecommand\transparent[1]{%
    \errmessage{(Inkscape) Transparency is used (non-zero) for the text in Inkscape, but the package 'transparent.sty' is not loaded}%
    \renewcommand\transparent[1]{}%
  }%
  \providecommand\rotatebox[2]{#2}%
  \newcommand*\fsize{\dimexpr\f@size pt\relax}%
  \newcommand*\lineheight[1]{\fontsize{\fsize}{#1\fsize}\selectfont}%
  \ifx\svgwidth\undefined%
    \setlength{\unitlength}{137.48101783bp}%
    \ifx\svgscale\undefined%
      \relax%
    \else%
      \setlength{\unitlength}{\unitlength * \real{\svgscale}}%
    \fi%
  \else%
    \setlength{\unitlength}{\svgwidth}%
  \fi%
  \global\let\svgwidth\undefined%
  \global\let\svgscale\undefined%
  \makeatother%
  \begin{picture}(1,1.00000181)%
    \lineheight{1}%
    \setlength\tabcolsep{0pt}%
    \put(0,0){\includegraphics[width=\unitlength,page=1]{TypeA_labels_s.pdf}}%
    \put(0.08086257,0.86601402){\color[rgb]{1,0,0}\makebox(0,0)[lt]{\lineheight{1.25}\smash{\begin{tabular}[t]{l}$\rho_1$\end{tabular}}}}%
    \put(0.08743294,0.10175109){\color[rgb]{1,0,0}\makebox(0,0)[lt]{\lineheight{1.25}\smash{\begin{tabular}[t]{l}$\rho_2$\end{tabular}}}}%
    \put(0.83781054,0.10175109){\color[rgb]{1,0,0}\makebox(0,0)[lt]{\lineheight{1.25}\smash{\begin{tabular}[t]{l}$\rho_3$\end{tabular}}}}%
    \put(0.44772905,0.79321015){\color[rgb]{1,0,0}\makebox(0,0)[lt]{\lineheight{1.25}\smash{\begin{tabular}[t]{l}$\rho_2\rho_1$\end{tabular}}}}%
    \put(0.12976368,0.34908476){\color[rgb]{1,0,0}\makebox(0,0)[lt]{\lineheight{1.25}\smash{\begin{tabular}[t]{l}$\rho_3\rho_2$\end{tabular}}}}%
    \put(0.76629872,0.86628761){\color[rgb]{1,0,0}\makebox(0,0)[lt]{\lineheight{1.25}\smash{\begin{tabular}[t]{l}$\rho_3\rho_2\rho_1$\end{tabular}}}}%
  \end{picture}%
\endgroup%

%% file: figures/Type_D_labels_s.pdf_tex
\begingroup%
  \makeatletter%
  \providecommand\color[2][]{%
    \errmessage{(Inkscape) Color is used for the text in Inkscape, but the package 'color.sty' is not loaded}%
    \renewcommand\color[2][]{}%
  }%
  \providecommand\transparent[1]{%
    \errmessage{(Inkscape) Transparency is used (non-zero) for the text in Inkscape, but the package 'transparent.sty' is not loaded}%
    \renewcommand\transparent[1]{}%
  }%
  \providecommand\rotatebox[2]{#2}%
  \newcommand*\fsize{\dimexpr\f@size pt\relax}%
  \newcommand*\lineheight[1]{\fontsize{\fsize}{#1\fsize}\selectfont}%
  \ifx\svgwidth\undefined%
    \setlength{\unitlength}{137.48100701bp}%
    \ifx\svgscale\undefined%
      \relax%
    \else%
      \setlength{\unitlength}{\unitlength * \real{\svgscale}}%
    \fi%
  \else%
    \setlength{\unitlength}{\svgwidth}%
  \fi%
  \global\let\svgwidth\undefined%
  \global\let\svgscale\undefined%
  \makeatother%
  \begin{picture}(1,1.00000157)%
    \lineheight{1}%
    \setlength\tabcolsep{0pt}%
    \put(0,0){\includegraphics[width=\unitlength,page=1]{Type_D_labels_s.pdf}}%
    \put(0.82888678,0.07998852){\color[rgb]{0,0,1}\makebox(0,0)[lt]{\lineheight{1.25}\smash{\begin{tabular}[t]{l}$\rho_1$\end{tabular}}}}%
    \put(0.86595822,0.87043701){\color[rgb]{0,0,1}\makebox(0,0)[lt]{\lineheight{1.25}\smash{\begin{tabular}[t]{l}$\rho_2$\end{tabular}}}}%
    \put(0.09812384,0.88789398){\color[rgb]{0,0,1}\makebox(0,0)[lt]{\lineheight{1.25}\smash{\begin{tabular}[t]{l}$\rho_3$\end{tabular}}}}%
    \put(0.40348848,0.87668254){\color[rgb]{0,0,1}\makebox(0,0)[lt]{\lineheight{1.25}\smash{\begin{tabular}[t]{l}$\rho_{12}$\end{tabular}}}}%
    \put(0.09973163,0.46090682){\color[rgb]{0,0,1}\makebox(0,0)[lt]{\lineheight{1.25}\smash{\begin{tabular}[t]{l}$\rho_{23}$\end{tabular}}}}%
    \put(0.07302618,0.08844342){\color[rgb]{0,0,1}\makebox(0,0)[lt]{\lineheight{1.25}\smash{\begin{tabular}[t]{l}$\rho_{123}$\end{tabular}}}}%
  \end{picture}%
\endgroup%

%% file: figures/Imm_AD_pairing.pdf_tex
\begingroup%
  \makeatletter%
  \providecommand\color[2][]{%
    \errmessage{(Inkscape) Color is used for the text in Inkscape, but the package 'color.sty' is not loaded}%
    \renewcommand\color[2][]{}%
  }%
  \providecommand\transparent[1]{%
    \errmessage{(Inkscape) Transparency is used (non-zero) for the text in Inkscape, but the package 'transparent.sty' is not loaded}%
    \renewcommand\transparent[1]{}%
  }%
  \providecommand\rotatebox[2]{#2}%
  \newcommand*\fsize{\dimexpr\f@size pt\relax}%
  \newcommand*\lineheight[1]{\fontsize{\fsize}{#1\fsize}\selectfont}%
  \ifx\svgwidth\undefined%
    \setlength{\unitlength}{169.79603685bp}%
    \ifx\svgscale\undefined%
      \relax%
    \else%
      \setlength{\unitlength}{\unitlength * \real{\svgscale}}%
    \fi%
  \else%
    \setlength{\unitlength}{\svgwidth}%
  \fi%
  \global\let\svgwidth\undefined%
  \global\let\svgscale\undefined%
  \makeatother%
  \begin{picture}(1,0.99999541)%
    \lineheight{1}%
    \setlength\tabcolsep{0pt}%
    \put(0,0){\includegraphics[width=\unitlength,page=1]{Imm_AD_pairing.pdf}}%
    \put(0.11592896,0.78340085){\color[rgb]{0,0,0}\makebox(0,0)[lt]{\lineheight{1.25}\smash{\begin{tabular}[t]{l}$b\boxtimes y$\end{tabular}}}}%
    \put(0.78420363,0.16988071){\color[rgb]{0,0,0}\makebox(0,0)[lt]{\lineheight{1.25}\smash{\begin{tabular}[t]{l}$a\boxtimes x$\end{tabular}}}}%
    \put(0.79842578,0.52266946){\color[rgb]{1,0,0}\makebox(0,0)[lt]{\lineheight{1.25}\smash{\begin{tabular}[t]{l}$a$\end{tabular}}}}%
    \put(0.43496546,0.17507618){\color[rgb]{0,0,1}\makebox(0,0)[lt]{\lineheight{1.25}\smash{\begin{tabular}[t]{l}$x$\end{tabular}}}}%
    \put(0.22044004,0.46280648){\color[rgb]{0,0,1}\makebox(0,0)[lt]{\lineheight{1.25}\smash{\begin{tabular}[t]{l}$y$\end{tabular}}}}%
    \put(0.52510141,0.7848916){\color[rgb]{1,0,0}\makebox(0,0)[lt]{\lineheight{1.25}\smash{\begin{tabular}[t]{l}$b$\end{tabular}}}}%
  \end{picture}%
\endgroup%

%% file: figures/Tpq_seifert.pdf_tex
\begingroup%
  \makeatletter%
  \providecommand\color[2][]{%
    \errmessage{(Inkscape) Color is used for the text in Inkscape, but the package 'color.sty' is not loaded}%
    \renewcommand\color[2][]{}%
  }%
  \providecommand\transparent[1]{%
    \errmessage{(Inkscape) Transparency is used (non-zero) for the text in Inkscape, but the package 'transparent.sty' is not loaded}%
    \renewcommand\transparent[1]{}%
  }%
  \providecommand\rotatebox[2]{#2}%
  \newcommand*\fsize{\dimexpr\f@size pt\relax}%
  \newcommand*\lineheight[1]{\fontsize{\fsize}{#1\fsize}\selectfont}%
  \ifx\svgwidth\undefined%
    \setlength{\unitlength}{81.19450126bp}%
    \ifx\svgscale\undefined%
      \relax%
    \else%
      \setlength{\unitlength}{\unitlength * \real{\svgscale}}%
    \fi%
  \else%
    \setlength{\unitlength}{\svgwidth}%
  \fi%
  \global\let\svgwidth\undefined%
  \global\let\svgscale\undefined%
  \makeatother%
  \begin{picture}(1,3.17801067)%
    \lineheight{1}%
    \setlength\tabcolsep{0pt}%
    \put(0,0){\includegraphics[width=\unitlength,page=1]{Tpq_seifert.pdf}}%
    \put(0.5336356,2.41331214){\color[rgb]{1,0,0}\makebox(0,0)[lt]{\lineheight{1.25}\smash{\begin{tabular}[t]{l}$c_A(\xi^{\textit{in}}_{p,q})$\end{tabular}}}}%
    \put(0.06706839,3.03044803){\color[rgb]{1,0,0}\makebox(0,0)[lt]{\lineheight{1.25}\smash{\begin{tabular}[t]{l}$c_A(\xi_{p,q}^{123})$\end{tabular}}}}%
    \put(0.50109464,1.62493462){\color[rgb]{0,0,0}\rotatebox{-90}{\makebox(0,0)[lt]{\lineheight{1.25}\smash{\begin{tabular}[t]{l}$\dots$\end{tabular}}}}}%
    \put(0.56658033,0.69407835){\color[rgb]{1,0,0}\makebox(0,0)[lt]{\lineheight{1.25}\smash{\begin{tabular}[t]{l}$\dots$\end{tabular}}}}%
    \put(0.08192555,0.6896892){\color[rgb]{1,0,0}\rotatebox{-90}{\makebox(0,0)[lt]{\lineheight{1.25}\smash{\begin{tabular}[t]{l}$\dots$\end{tabular}}}}}%
  \end{picture}%
\endgroup%

%% file: figures/LHT_typeA_beta.pdf_tex
\begingroup%
  \makeatletter%
  \providecommand\color[2][]{%
    \errmessage{(Inkscape) Color is used for the text in Inkscape, but the package 'color.sty' is not loaded}%
    \renewcommand\color[2][]{}%
  }%
  \providecommand\transparent[1]{%
    \errmessage{(Inkscape) Transparency is used (non-zero) for the text in Inkscape, but the package 'transparent.sty' is not loaded}%
    \renewcommand\transparent[1]{}%
  }%
  \providecommand\rotatebox[2]{#2}%
  \newcommand*\fsize{\dimexpr\f@size pt\relax}%
  \newcommand*\lineheight[1]{\fontsize{\fsize}{#1\fsize}\selectfont}%
  \ifx\svgwidth\undefined%
    \setlength{\unitlength}{120.68121362bp}%
    \ifx\svgscale\undefined%
      \relax%
    \else%
      \setlength{\unitlength}{\unitlength * \real{\svgscale}}%
    \fi%
  \else%
    \setlength{\unitlength}{\svgwidth}%
  \fi%
  \global\let\svgwidth\undefined%
  \global\let\svgscale\undefined%
  \makeatother%
  \begin{picture}(1,1.0485433)%
    \lineheight{1}%
    \setlength\tabcolsep{0pt}%
    \put(0,0){\includegraphics[width=\unitlength,page=1]{LHT_typeA_beta.pdf}}%
    \put(0.22779001,0.11680895){\color[rgb]{0,0,0}\makebox(0,0)[lt]{\lineheight{1.25}\smash{\begin{tabular}[t]{l}1\end{tabular}}}}%
    \put(0.75293336,0.11912871){\color[rgb]{0,0,0}\makebox(0,0)[lt]{\lineheight{1.25}\smash{\begin{tabular}[t]{l}3\end{tabular}}}}%
    \put(0.90711069,0.51593375){\color[rgb]{0,0,0}\makebox(0,0)[lt]{\lineheight{1.25}\smash{\begin{tabular}[t]{l}321\end{tabular}}}}%
    \put(0.04121034,0.49822647){\color[rgb]{0,0,0}\makebox(0,0)[lt]{\lineheight{1.25}\smash{\begin{tabular}[t]{l}2\end{tabular}}}}%
    \put(0.25753425,0.88260215){\color[rgb]{0,0,0}\makebox(0,0)[lt]{\lineheight{1.25}\smash{\begin{tabular}[t]{l}3\end{tabular}}}}%
    \put(0.7089164,0.87706359){\color[rgb]{0,0,0}\makebox(0,0)[lt]{\lineheight{1.25}\smash{\begin{tabular}[t]{l}1\end{tabular}}}}%
    \put(0.93032643,0.24054199){\color[rgb]{0,0,0}\makebox(0,0)[lt]{\lineheight{1.25}\smash{\begin{tabular}[t]{l}$T$\end{tabular}}}}%
    \put(0.010548,0.2281125){\color[rgb]{0,0,0}\makebox(0,0)[lt]{\lineheight{1.25}\smash{\begin{tabular}[t]{l}$B$\end{tabular}}}}%
    \put(0.43172787,1.02941403){\color[rgb]{0,0,0}\makebox(0,0)[lt]{\lineheight{1.25}\smash{\begin{tabular}[t]{l}$M=c_A(\xi^{\textit{in}}_{-2,3})$\end{tabular}}}}%
    \put(0.49191496,0.0072636){\color[rgb]{0,0,0}\makebox(0,0)[lt]{\lineheight{1.25}\smash{\begin{tabular}[t]{l}32\end{tabular}}}}%
    \put(-0.00188149,0.75014873){\color[rgb]{0,0,0}\makebox(0,0)[lt]{\lineheight{1.25}\smash{\begin{tabular}[t]{l}$d$\end{tabular}}}}%
    \put(0.95518598,0.75014873){\color[rgb]{0,0,0}\makebox(0,0)[lt]{\lineheight{1.25}\smash{\begin{tabular}[t]{l}$a$\end{tabular}}}}%
    \put(0.65687969,0.00438214){\color[rgb]{0,0,0}\makebox(0,0)[lt]{\lineheight{1.25}\smash{\begin{tabular}[t]{l}$b$\end{tabular}}}}%
    \put(0.29642581,0.00438214){\color[rgb]{0,0,0}\makebox(0,0)[lt]{\lineheight{1.25}\smash{\begin{tabular}[t]{l}$c$\end{tabular}}}}%
  \end{picture}%
\endgroup%

%% file: figures/LHT_typeD.pdf_tex
\begingroup%
  \makeatletter%
  \providecommand\color[2][]{%
    \errmessage{(Inkscape) Color is used for the text in Inkscape, but the package 'color.sty' is not loaded}%
    \renewcommand\color[2][]{}%
  }%
  \providecommand\transparent[1]{%
    \errmessage{(Inkscape) Transparency is used (non-zero) for the text in Inkscape, but the package 'transparent.sty' is not loaded}%
    \renewcommand\transparent[1]{}%
  }%
  \providecommand\rotatebox[2]{#2}%
  \newcommand*\fsize{\dimexpr\f@size pt\relax}%
  \newcommand*\lineheight[1]{\fontsize{\fsize}{#1\fsize}\selectfont}%
  \ifx\svgwidth\undefined%
    \setlength{\unitlength}{134.24129456bp}%
    \ifx\svgscale\undefined%
      \relax%
    \else%
      \setlength{\unitlength}{\unitlength * \real{\svgscale}}%
    \fi%
  \else%
    \setlength{\unitlength}{\svgwidth}%
  \fi%
  \global\let\svgwidth\undefined%
  \global\let\svgscale\undefined%
  \makeatother%
  \begin{picture}(1,0.87434886)%
    \lineheight{1}%
    \setlength\tabcolsep{0pt}%
    \put(0,0){\includegraphics[width=\unitlength,page=1]{LHT_typeD.pdf}}%
    \put(0.1653903,0.098931){\color[rgb]{0,0,0}\makebox(0,0)[lt]{\lineheight{1.25}\smash{\begin{tabular}[t]{l}1\end{tabular}}}}%
    \put(0.63748747,0.10101644){\color[rgb]{0,0,0}\makebox(0,0)[lt]{\lineheight{1.25}\smash{\begin{tabular}[t]{l}3\end{tabular}}}}%
    \put(0.77609092,0.45773911){\color[rgb]{0,0,0}\makebox(0,0)[lt]{\lineheight{1.25}\smash{\begin{tabular}[t]{l}123\end{tabular}}}}%
    \put(-0.00234245,0.44182049){\color[rgb]{0,0,0}\makebox(0,0)[lt]{\lineheight{1.25}\smash{\begin{tabular}[t]{l}2\end{tabular}}}}%
    \put(0.19212999,0.78736933){\color[rgb]{0,0,0}\makebox(0,0)[lt]{\lineheight{1.25}\smash{\begin{tabular}[t]{l}3\end{tabular}}}}%
    \put(0.59791679,0.78239023){\color[rgb]{0,0,0}\makebox(0,0)[lt]{\lineheight{1.25}\smash{\begin{tabular}[t]{l}1\end{tabular}}}}%
    \put(0.75096613,0.69301395){\color[rgb]{0,0,0}\makebox(0,0)[lt]{\lineheight{1.25}\smash{\begin{tabular}[t]{l}$c_D(\xi^{\textit{in}}_{-2,3})$\end{tabular}}}}%
    \put(0.40283526,0.00045114){\color[rgb]{0,0,0}\makebox(0,0)[lt]{\lineheight{1.25}\smash{\begin{tabular}[t]{l}23\end{tabular}}}}%
  \end{picture}%
\endgroup%

%% file: figures/LHT_D_full.pdf_tex
\begingroup%
  \makeatletter%
  \providecommand\color[2][]{%
    \errmessage{(Inkscape) Color is used for the text in Inkscape, but the package 'color.sty' is not loaded}%
    \renewcommand\color[2][]{}%
  }%
  \providecommand\transparent[1]{%
    \errmessage{(Inkscape) Transparency is used (non-zero) for the text in Inkscape, but the package 'transparent.sty' is not loaded}%
    \renewcommand\transparent[1]{}%
  }%
  \providecommand\rotatebox[2]{#2}%
  \newcommand*\fsize{\dimexpr\f@size pt\relax}%
  \newcommand*\lineheight[1]{\fontsize{\fsize}{#1\fsize}\selectfont}%
  \ifx\svgwidth\undefined%
    \setlength{\unitlength}{206.90863181bp}%
    \ifx\svgscale\undefined%
      \relax%
    \else%
      \setlength{\unitlength}{\unitlength * \real{\svgscale}}%
    \fi%
  \else%
    \setlength{\unitlength}{\svgwidth}%
  \fi%
  \global\let\svgwidth\undefined%
  \global\let\svgscale\undefined%
  \makeatother%
  \begin{picture}(1,1.51628274)%
    \lineheight{1}%
    \setlength\tabcolsep{0pt}%
    \put(0,0){\includegraphics[width=\unitlength,page=1]{LHT_D_full.pdf}}%
    \put(0.08219959,0.69006016){\color[rgb]{0,0,0}\makebox(0,0)[lt]{\lineheight{1.25}\smash{\begin{tabular}[t]{l}$\iota_0^{\vee}T$\end{tabular}}}}%
    \put(0.16289882,1.07196861){\color[rgb]{0,0,0}\makebox(0,0)[lt]{\lineheight{1.25}\smash{\begin{tabular}[t]{l}$\iota_1^{\vee}c$\end{tabular}}}}%
    \put(0.23461333,0.69006016){\color[rgb]{0,0,0}\makebox(0,0)[lt]{\lineheight{1.25}\smash{\begin{tabular}[t]{l}$\rho_1^{\vee}T$\end{tabular}}}}%
    \put(0.39934711,0.69006016){\color[rgb]{0,0,0}\makebox(0,0)[lt]{\lineheight{1.25}\smash{\begin{tabular}[t]{l}$\rho_3^{\vee}T$\end{tabular}}}}%
    \put(0.57245449,0.69006016){\color[rgb]{0,0,0}\makebox(0,0)[lt]{\lineheight{1.25}\smash{\begin{tabular}[t]{l}$\rho_{12}^{\vee}T$\end{tabular}}}}%
    \put(0.73218527,0.69006016){\color[rgb]{0,0,0}\makebox(0,0)[lt]{\lineheight{1.25}\smash{\begin{tabular}[t]{l}$\rho_{123}^{\vee}T$\end{tabular}}}}%
    \put(0.45163802,1.07196861){\color[rgb]{0,0,0}\makebox(0,0)[lt]{\lineheight{1.25}\smash{\begin{tabular}[t]{l}$\rho_2^{\vee}c$\end{tabular}}}}%
    \put(0.7184219,1.06036928){\color[rgb]{0,0,0}\makebox(0,0)[lt]{\lineheight{1.25}\smash{\begin{tabular}[t]{l}$\rho_{23}^{\vee}c$\end{tabular}}}}%
    \put(0.61555669,1.03594362){\color[rgb]{0,0,0}\makebox(0,0)[lt]{\lineheight{1.25}\smash{\begin{tabular}[t]{l}3\end{tabular}}}}%
    \put(0.34156237,1.03594362){\color[rgb]{0,0,0}\makebox(0,0)[lt]{\lineheight{1.25}\smash{\begin{tabular}[t]{l}2\end{tabular}}}}%
    \put(0.69113941,0.66604956){\color[rgb]{0,0,0}\makebox(0,0)[lt]{\lineheight{1.25}\smash{\begin{tabular}[t]{l}3\end{tabular}}}}%
    \put(0.19250886,0.66624707){\color[rgb]{0,0,0}\makebox(0,0)[lt]{\lineheight{1.25}\smash{\begin{tabular}[t]{l}1\end{tabular}}}}%
    \put(0.22219462,0.54828812){\color[rgb]{0,0,0}\makebox(0,0)[lt]{\lineheight{1.25}\smash{\begin{tabular}[t]{l}3\end{tabular}}}}%
    \put(0.48403998,0.56129717){\color[rgb]{0,0,0}\makebox(0,0)[lt]{\lineheight{1.25}\smash{\begin{tabular}[t]{l}2\end{tabular}}}}%
    \put(0.08219959,0.41606451){\color[rgb]{0,0,0}\makebox(0,0)[lt]{\lineheight{1.25}\smash{\begin{tabular}[t]{l}$\iota_0^{\vee}M$\end{tabular}}}}%
    \put(0.23461333,0.41606451){\color[rgb]{0,0,0}\makebox(0,0)[lt]{\lineheight{1.25}\smash{\begin{tabular}[t]{l}$\rho_1^{\vee}M$\end{tabular}}}}%
    \put(0.39934711,0.41606451){\color[rgb]{0,0,0}\makebox(0,0)[lt]{\lineheight{1.25}\smash{\begin{tabular}[t]{l}$\rho_3^{\vee}M$\end{tabular}}}}%
    \put(0.57245449,0.41606451){\color[rgb]{0,0,0}\makebox(0,0)[lt]{\lineheight{1.25}\smash{\begin{tabular}[t]{l}$\rho_{12}^{\vee}M$\end{tabular}}}}%
    \put(0.73218527,0.41606451){\color[rgb]{0,0,0}\makebox(0,0)[lt]{\lineheight{1.25}\smash{\begin{tabular}[t]{l}$\rho_{123}^{\vee}M$\end{tabular}}}}%
    \put(0.69113941,0.39205391){\color[rgb]{0,0,0}\makebox(0,0)[lt]{\lineheight{1.25}\smash{\begin{tabular}[t]{l}3\end{tabular}}}}%
    \put(0.19250886,0.39225142){\color[rgb]{0,0,0}\makebox(0,0)[lt]{\lineheight{1.25}\smash{\begin{tabular}[t]{l}1\end{tabular}}}}%
    \put(0.22219462,0.27429247){\color[rgb]{0,0,0}\makebox(0,0)[lt]{\lineheight{1.25}\smash{\begin{tabular}[t]{l}3\end{tabular}}}}%
    \put(0.48403998,0.2873016){\color[rgb]{0,0,0}\makebox(0,0)[lt]{\lineheight{1.25}\smash{\begin{tabular}[t]{l}2\end{tabular}}}}%
    \put(0.08219959,0.14206468){\color[rgb]{0,0,0}\makebox(0,0)[lt]{\lineheight{1.25}\smash{\begin{tabular}[t]{l}$\iota_0^{\vee}B$\end{tabular}}}}%
    \put(0.23461333,0.14206468){\color[rgb]{0,0,0}\makebox(0,0)[lt]{\lineheight{1.25}\smash{\begin{tabular}[t]{l}$\rho_1^{\vee}B$\end{tabular}}}}%
    \put(0.39934711,0.14206468){\color[rgb]{0,0,0}\makebox(0,0)[lt]{\lineheight{1.25}\smash{\begin{tabular}[t]{l}$\rho_3^{\vee}B$\end{tabular}}}}%
    \put(0.57245449,0.14206468){\color[rgb]{0,0,0}\makebox(0,0)[lt]{\lineheight{1.25}\smash{\begin{tabular}[t]{l}$\rho_{12}^{\vee}B$\end{tabular}}}}%
    \put(0.73218527,0.14206468){\color[rgb]{0,0,0}\makebox(0,0)[lt]{\lineheight{1.25}\smash{\begin{tabular}[t]{l}$\rho_{123}^{\vee}B$\end{tabular}}}}%
    \put(0.69113941,0.11805408){\color[rgb]{0,0,0}\makebox(0,0)[lt]{\lineheight{1.25}\smash{\begin{tabular}[t]{l}3\end{tabular}}}}%
    \put(0.19250886,0.11825158){\color[rgb]{0,0,0}\makebox(0,0)[lt]{\lineheight{1.25}\smash{\begin{tabular}[t]{l}1\end{tabular}}}}%
    \put(0.22219462,0.00029272){\color[rgb]{0,0,0}\makebox(0,0)[lt]{\lineheight{1.25}\smash{\begin{tabular}[t]{l}3\end{tabular}}}}%
    \put(0.48403998,0.01330185){\color[rgb]{0,0,0}\makebox(0,0)[lt]{\lineheight{1.25}\smash{\begin{tabular}[t]{l}2\end{tabular}}}}%
    \put(0.16289882,0.8938688){\color[rgb]{0,0,0}\makebox(0,0)[lt]{\lineheight{1.25}\smash{\begin{tabular}[t]{l}$\iota_1^{\vee}d$\end{tabular}}}}%
    \put(0.45163802,0.8938688){\color[rgb]{0,0,0}\makebox(0,0)[lt]{\lineheight{1.25}\smash{\begin{tabular}[t]{l}$\rho_2^{\vee}d$\end{tabular}}}}%
    \put(0.7184219,0.88226946){\color[rgb]{0,0,0}\makebox(0,0)[lt]{\lineheight{1.25}\smash{\begin{tabular}[t]{l}$\rho_{23}^{\vee}d$\end{tabular}}}}%
    \put(0.61555669,0.85784378){\color[rgb]{0,0,0}\makebox(0,0)[lt]{\lineheight{1.25}\smash{\begin{tabular}[t]{l}3\end{tabular}}}}%
    \put(0.34156237,0.85784378){\color[rgb]{0,0,0}\makebox(0,0)[lt]{\lineheight{1.25}\smash{\begin{tabular}[t]{l}2\end{tabular}}}}%
    \put(0.16289924,1.42816824){\color[rgb]{0,0,0}\makebox(0,0)[lt]{\lineheight{1.25}\smash{\begin{tabular}[t]{l}$\iota_1^{\vee}a$\end{tabular}}}}%
    \put(0.45163827,1.42816824){\color[rgb]{0,0,0}\makebox(0,0)[lt]{\lineheight{1.25}\smash{\begin{tabular}[t]{l}$\rho_2^{\vee}a$\end{tabular}}}}%
    \put(0.71842207,1.4165689){\color[rgb]{0,0,0}\makebox(0,0)[lt]{\lineheight{1.25}\smash{\begin{tabular}[t]{l}$\rho_{23}^{\vee}a$\end{tabular}}}}%
    \put(0.61555685,1.39214324){\color[rgb]{0,0,0}\makebox(0,0)[lt]{\lineheight{1.25}\smash{\begin{tabular}[t]{l}3\end{tabular}}}}%
    \put(0.3415627,1.39214324){\color[rgb]{0,0,0}\makebox(0,0)[lt]{\lineheight{1.25}\smash{\begin{tabular}[t]{l}2\end{tabular}}}}%
    \put(0.16289924,1.2500685){\color[rgb]{0,0,0}\makebox(0,0)[lt]{\lineheight{1.25}\smash{\begin{tabular}[t]{l}$\iota_1^{\vee}b$\end{tabular}}}}%
    \put(0.45163827,1.2500685){\color[rgb]{0,0,0}\makebox(0,0)[lt]{\lineheight{1.25}\smash{\begin{tabular}[t]{l}$\rho_2^{\vee}b$\end{tabular}}}}%
    \put(0.71842207,1.23846918){\color[rgb]{0,0,0}\makebox(0,0)[lt]{\lineheight{1.25}\smash{\begin{tabular}[t]{l}$\rho_{23}^{\vee}b$\end{tabular}}}}%
    \put(0.61555685,1.21404352){\color[rgb]{0,0,0}\makebox(0,0)[lt]{\lineheight{1.25}\smash{\begin{tabular}[t]{l}3\end{tabular}}}}%
    \put(0.3415627,1.21404352){\color[rgb]{0,0,0}\makebox(0,0)[lt]{\lineheight{1.25}\smash{\begin{tabular}[t]{l}2\end{tabular}}}}%
  \end{picture}%
\endgroup%

%% file: figures/Trefoil_splice.pdf_tex
\begingroup%
  \makeatletter%
  \providecommand\color[2][]{%
    \errmessage{(Inkscape) Color is used for the text in Inkscape, but the package 'color.sty' is not loaded}%
    \renewcommand\color[2][]{}%
  }%
  \providecommand\transparent[1]{%
    \errmessage{(Inkscape) Transparency is used (non-zero) for the text in Inkscape, but the package 'transparent.sty' is not loaded}%
    \renewcommand\transparent[1]{}%
  }%
  \providecommand\rotatebox[2]{#2}%
  \newcommand*\fsize{\dimexpr\f@size pt\relax}%
  \newcommand*\lineheight[1]{\fontsize{\fsize}{#1\fsize}\selectfont}%
  \ifx\svgwidth\undefined%
    \setlength{\unitlength}{232.42094025bp}%
    \ifx\svgscale\undefined%
      \relax%
    \else%
      \setlength{\unitlength}{\unitlength * \real{\svgscale}}%
    \fi%
  \else%
    \setlength{\unitlength}{\svgwidth}%
  \fi%
  \global\let\svgwidth\undefined%
  \global\let\svgscale\undefined%
  \makeatother%
  \begin{picture}(1,0.88411934)%
    \lineheight{1}%
    \setlength\tabcolsep{0pt}%
    \put(0,0){\includegraphics[width=\unitlength,page=1]{Trefoil_splice.pdf}}%
    \put(0.72146067,0.86977711){\color[rgb]{1,0,0}\makebox(0,0)[lt]{\lineheight{1.25}\smash{\begin{tabular}[t]{l}$c_A(\xi_{p,q}^{123})$\end{tabular}}}}%
    \put(0.7189287,0.05305885){\color[rgb]{0.10196078,0.10196078,0.10196078}\makebox(0,0)[lt]{\lineheight{1.25}\smash{\begin{tabular}[t]{l}$c(\xi)$\end{tabular}}}}%
    \put(0.22824648,0.65779865){\color[rgb]{0.10196078,0.10196078,0.10196078}\makebox(0,0)[lt]{\lineheight{1.25}\smash{\begin{tabular}[t]{l}$x$\end{tabular}}}}%
    \put(0.81860262,0.18579196){\color[rgb]{0.10196078,0.10196078,0.10196078}\makebox(0,0)[lt]{\lineheight{1.25}\smash{\begin{tabular}[t]{l}$y$\end{tabular}}}}%
    \put(-0.00141056,0.11064191){\color[rgb]{0,0,1}\makebox(0,0)[lt]{\lineheight{1.25}\smash{\begin{tabular}[t]{l}$c_D(\xi^{\textit{in}}_{-2,3})$\end{tabular}}}}%
    \put(0.87819394,0.66961155){\color[rgb]{1,0,0}\makebox(0,0)[lt]{\lineheight{1.25}\smash{\begin{tabular}[t]{l}$\dots$\end{tabular}}}}%
    \put(0.69404674,0.6987308){\color[rgb]{1,0,0}\makebox(0,0)[lt]{\lineheight{1.25}\smash{\begin{tabular}[t]{l}$\dots$\end{tabular}}}}%
    \put(0.71009136,0.78468613){\color[rgb]{1,0,0}\rotatebox{-90}{\makebox(0,0)[lt]{\lineheight{1.25}\smash{\begin{tabular}[t]{l}$\dots$\end{tabular}}}}}%
  \end{picture}%
\endgroup%

%% file: figures/Fig8_typeA_beta.pdf_tex
\begingroup%
  \makeatletter%
  \providecommand\color[2][]{%
    \errmessage{(Inkscape) Color is used for the text in Inkscape, but the package 'color.sty' is not loaded}%
    \renewcommand\color[2][]{}%
  }%
  \providecommand\transparent[1]{%
    \errmessage{(Inkscape) Transparency is used (non-zero) for the text in Inkscape, but the package 'transparent.sty' is not loaded}%
    \renewcommand\transparent[1]{}%
  }%
  \providecommand\rotatebox[2]{#2}%
  \newcommand*\fsize{\dimexpr\f@size pt\relax}%
  \newcommand*\lineheight[1]{\fontsize{\fsize}{#1\fsize}\selectfont}%
  \ifx\svgwidth\undefined%
    \setlength{\unitlength}{144.06606719bp}%
    \ifx\svgscale\undefined%
      \relax%
    \else%
      \setlength{\unitlength}{\unitlength * \real{\svgscale}}%
    \fi%
  \else%
    \setlength{\unitlength}{\svgwidth}%
  \fi%
  \global\let\svgwidth\undefined%
  \global\let\svgscale\undefined%
  \makeatother%
  \begin{picture}(1,0.95458644)%
    \lineheight{1}%
    \setlength\tabcolsep{0pt}%
    \put(0,0){\includegraphics[width=\unitlength,page=1]{Fig8_typeA_beta.pdf}}%
    \put(0.21163898,0.11942871){\color[rgb]{0,0,0}\makebox(0,0)[lt]{\lineheight{1.25}\smash{\begin{tabular}[t]{l}1\end{tabular}}}}%
    \put(0.94928515,0.48734047){\color[rgb]{0,0,0}\makebox(0,0)[lt]{\lineheight{1.25}\smash{\begin{tabular}[t]{l}3\end{tabular}}}}%
    \put(0.85357548,0.7869481){\color[rgb]{0,0,0}\makebox(0,0)[lt]{\lineheight{1.25}\smash{\begin{tabular}[t]{l}321\end{tabular}}}}%
    \put(0.03452113,0.43893441){\color[rgb]{0,0,0}\makebox(0,0)[lt]{\lineheight{1.25}\smash{\begin{tabular}[t]{l}2\end{tabular}}}}%
    \put(0.19763772,0.80557595){\color[rgb]{0,0,0}\makebox(0,0)[lt]{\lineheight{1.25}\smash{\begin{tabular}[t]{l}3\end{tabular}}}}%
    \put(0.51054972,0.92286871){\color[rgb]{0,0,0}\makebox(0,0)[lt]{\lineheight{1.25}\smash{\begin{tabular}[t]{l}1\end{tabular}}}}%
    \put(0.96213485,0.64032903){\color[rgb]{0,0,0}\makebox(0,0)[lt]{\lineheight{1.25}\smash{\begin{tabular}[t]{l}$T$\end{tabular}}}}%
    \put(0.00883591,0.24390122){\color[rgb]{0,0,0}\makebox(0,0)[lt]{\lineheight{1.25}\smash{\begin{tabular}[t]{l}$B$\end{tabular}}}}%
    \put(0.31479622,0.93856184){\color[rgb]{0,0,0}\makebox(0,0)[lt]{\lineheight{1.25}\smash{\begin{tabular}[t]{l}$E$\end{tabular}}}}%
    \put(0.82333706,0.12240413){\color[rgb]{0,0,0}\makebox(0,0)[lt]{\lineheight{1.25}\smash{\begin{tabular}[t]{l}2\end{tabular}}}}%
    \put(-0.00157602,0.64996456){\color[rgb]{0,0,0}\makebox(0,0)[lt]{\lineheight{1.25}\smash{\begin{tabular}[t]{l}$d$\end{tabular}}}}%
    \put(0.26830515,0.02525106){\color[rgb]{0,0,0}\makebox(0,0)[lt]{\lineheight{1.25}\smash{\begin{tabular}[t]{l}$a$\end{tabular}}}}%
    \put(0.69578813,0.93108594){\color[rgb]{0,0,0}\makebox(0,0)[lt]{\lineheight{1.25}\smash{\begin{tabular}[t]{l}$b$\end{tabular}}}}%
    \put(0.95130992,0.25752443){\color[rgb]{0,0,0}\makebox(0,0)[lt]{\lineheight{1.25}\smash{\begin{tabular}[t]{l}$c$\end{tabular}}}}%
    \put(0.68697331,0.00367049){\color[rgb]{0,0,0}\makebox(0,0)[lt]{\lineheight{1.25}\smash{\begin{tabular}[t]{l}$I$\end{tabular}}}}%
    \put(0.45657135,0.56899484){\color[rgb]{0,0,0}\makebox(0,0)[lt]{\lineheight{1.25}\smash{\begin{tabular}[t]{l}$Z$\end{tabular}}}}%
    \put(0.47749488,0.01850208){\color[rgb]{0,0,0}\makebox(0,0)[lt]{\lineheight{1.25}\smash{\begin{tabular}[t]{l}321\end{tabular}}}}%
    \put(0.51358338,0.33741848){\color[rgb]{0,0,0}\makebox(0,0)[lt]{\lineheight{1.25}\smash{\begin{tabular}[t]{l}21\end{tabular}}}}%
  \end{picture}%
\endgroup%

%% file: figures/Fig8_typeD.pdf_tex
\begingroup%
  \makeatletter%
  \providecommand\color[2][]{%
    \errmessage{(Inkscape) Color is used for the text in Inkscape, but the package 'color.sty' is not loaded}%
    \renewcommand\color[2][]{}%
  }%
  \providecommand\transparent[1]{%
    \errmessage{(Inkscape) Transparency is used (non-zero) for the text in Inkscape, but the package 'transparent.sty' is not loaded}%
    \renewcommand\transparent[1]{}%
  }%
  \providecommand\rotatebox[2]{#2}%
  \newcommand*\fsize{\dimexpr\f@size pt\relax}%
  \newcommand*\lineheight[1]{\fontsize{\fsize}{#1\fsize}\selectfont}%
  \ifx\svgwidth\undefined%
    \setlength{\unitlength}{133.7141887bp}%
    \ifx\svgscale\undefined%
      \relax%
    \else%
      \setlength{\unitlength}{\unitlength * \real{\svgscale}}%
    \fi%
  \else%
    \setlength{\unitlength}{\svgwidth}%
  \fi%
  \global\let\svgwidth\undefined%
  \global\let\svgscale\undefined%
  \makeatother%
  \begin{picture}(1,1.00095207)%
    \lineheight{1}%
    \setlength\tabcolsep{0pt}%
    \put(0,0){\includegraphics[width=\unitlength,page=1]{Fig8_typeD.pdf}}%
    \put(0.15482429,0.12041118){\color[rgb]{0,0,0}\makebox(0,0)[lt]{\lineheight{1.25}\smash{\begin{tabular}[t]{l}1\end{tabular}}}}%
    \put(0.98323139,0.50558792){\color[rgb]{0,0,0}\makebox(0,0)[lt]{\lineheight{1.25}\smash{\begin{tabular}[t]{l}3\end{tabular}}}}%
    \put(0.88011206,0.82839056){\color[rgb]{0,0,0}\makebox(0,0)[lt]{\lineheight{1.25}\smash{\begin{tabular}[t]{l}123\end{tabular}}}}%
    \put(-0.0023518,0.48708835){\color[rgb]{0,0,0}\makebox(0,0)[lt]{\lineheight{1.25}\smash{\begin{tabular}[t]{l}2\end{tabular}}}}%
    \put(0.17339295,0.84846053){\color[rgb]{0,0,0}\makebox(0,0)[lt]{\lineheight{1.25}\smash{\begin{tabular}[t]{l}3\end{tabular}}}}%
    \put(0.51052995,0.97483386){\color[rgb]{0,0,0}\makebox(0,0)[lt]{\lineheight{1.25}\smash{\begin{tabular}[t]{l}1\end{tabular}}}}%
    \put(0.84753265,0.11239895){\color[rgb]{0,0,0}\makebox(0,0)[lt]{\lineheight{1.25}\smash{\begin{tabular}[t]{l}2\end{tabular}}}}%
    \put(0.71010914,0.98368725){\color[rgb]{0,0,0}\makebox(0,0)[lt]{\lineheight{1.25}\smash{\begin{tabular}[t]{l}$\iota_0^{\vee}E$\end{tabular}}}}%
    \put(0.45237268,0.5935638){\color[rgb]{0,0,0}\makebox(0,0)[lt]{\lineheight{1.25}\smash{\begin{tabular}[t]{l}$Z$\end{tabular}}}}%
    \put(0.47491608,0.00045302){\color[rgb]{0,0,0}\makebox(0,0)[lt]{\lineheight{1.25}\smash{\begin{tabular}[t]{l}123\end{tabular}}}}%
    \put(0.51379848,0.34405927){\color[rgb]{0,0,0}\makebox(0,0)[lt]{\lineheight{1.25}\smash{\begin{tabular}[t]{l}12\end{tabular}}}}%
  \end{picture}%
\endgroup%

%% file: figures/FIg8_D_full.pdf_tex
\begingroup%
  \makeatletter%
  \providecommand\color[2][]{%
    \errmessage{(Inkscape) Color is used for the text in Inkscape, but the package 'color.sty' is not loaded}%
    \renewcommand\color[2][]{}%
  }%
  \providecommand\transparent[1]{%
    \errmessage{(Inkscape) Transparency is used (non-zero) for the text in Inkscape, but the package 'transparent.sty' is not loaded}%
    \renewcommand\transparent[1]{}%
  }%
  \providecommand\rotatebox[2]{#2}%
  \newcommand*\fsize{\dimexpr\f@size pt\relax}%
  \newcommand*\lineheight[1]{\fontsize{\fsize}{#1\fsize}\selectfont}%
  \ifx\svgwidth\undefined%
    \setlength{\unitlength}{222.33629122bp}%
    \ifx\svgscale\undefined%
      \relax%
    \else%
      \setlength{\unitlength}{\unitlength * \real{\svgscale}}%
    \fi%
  \else%
    \setlength{\unitlength}{\svgwidth}%
  \fi%
  \global\let\svgwidth\undefined%
  \global\let\svgscale\undefined%
  \makeatother%
  \begin{picture}(1,1.65293086)%
    \lineheight{1}%
    \setlength\tabcolsep{0pt}%
    \put(0,0){\includegraphics[width=\unitlength,page=1]{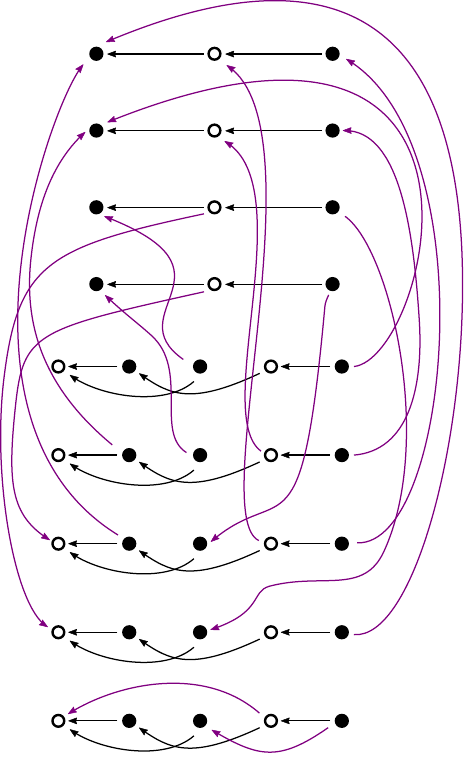}}%
    \put(0.10061551,0.89716053){\color[rgb]{0,0,0}\makebox(0,0)[lt]{\lineheight{1.25}\smash{\begin{tabular}[t]{l}$\iota_0^{\vee}T$\end{tabular}}}}%
    \put(0.17571512,1.2525688){\color[rgb]{0,0,0}\makebox(0,0)[lt]{\lineheight{1.25}\smash{\begin{tabular}[t]{l}$\iota_1^{\vee}c$\end{tabular}}}}%
    \put(0.24245344,0.89716053){\color[rgb]{0,0,0}\makebox(0,0)[lt]{\lineheight{1.25}\smash{\begin{tabular}[t]{l}$\rho_1^{\vee}T$\end{tabular}}}}%
    \put(0.39575653,0.89716053){\color[rgb]{0,0,0}\makebox(0,0)[lt]{\lineheight{1.25}\smash{\begin{tabular}[t]{l}$\rho_3^{\vee}T$\end{tabular}}}}%
    \put(0.55685219,0.89716053){\color[rgb]{0,0,0}\makebox(0,0)[lt]{\lineheight{1.25}\smash{\begin{tabular}[t]{l}$\rho_{12}^{\vee}T$\end{tabular}}}}%
    \put(0.70549943,0.89716053){\color[rgb]{0,0,0}\makebox(0,0)[lt]{\lineheight{1.25}\smash{\begin{tabular}[t]{l}$\rho_{123}^{\vee}T$\end{tabular}}}}%
    \put(0.44441903,1.2525688){\color[rgb]{0,0,0}\makebox(0,0)[lt]{\lineheight{1.25}\smash{\begin{tabular}[t]{l}$\rho_2^{\vee}c$\end{tabular}}}}%
    \put(0.69269109,1.24177433){\color[rgb]{0,0,0}\makebox(0,0)[lt]{\lineheight{1.25}\smash{\begin{tabular}[t]{l}$\rho_{23}^{\vee}c$\end{tabular}}}}%
    \put(0.59696357,1.21904354){\color[rgb]{0,0,0}\makebox(0,0)[lt]{\lineheight{1.25}\smash{\begin{tabular}[t]{l}3\end{tabular}}}}%
    \put(0.34198141,1.21904354){\color[rgb]{0,0,0}\makebox(0,0)[lt]{\lineheight{1.25}\smash{\begin{tabular}[t]{l}2\end{tabular}}}}%
    \put(0.6673017,0.874816){\color[rgb]{0,0,0}\makebox(0,0)[lt]{\lineheight{1.25}\smash{\begin{tabular}[t]{l}3\end{tabular}}}}%
    \put(0.20327055,0.8749998){\color[rgb]{0,0,0}\makebox(0,0)[lt]{\lineheight{1.25}\smash{\begin{tabular}[t]{l}1\end{tabular}}}}%
    \put(0.23089645,0.76522589){\color[rgb]{0,0,0}\makebox(0,0)[lt]{\lineheight{1.25}\smash{\begin{tabular}[t]{l}3\end{tabular}}}}%
    \put(0.47457266,0.77733225){\color[rgb]{0,0,0}\makebox(0,0)[lt]{\lineheight{1.25}\smash{\begin{tabular}[t]{l}2\end{tabular}}}}%
    \put(0.10061551,0.70592385){\color[rgb]{0,0,0}\makebox(0,0)[lt]{\lineheight{1.25}\smash{\begin{tabular}[t]{l}$\iota_0^{\vee}E$\end{tabular}}}}%
    \put(0.24245344,0.70592385){\color[rgb]{0,0,0}\makebox(0,0)[lt]{\lineheight{1.25}\smash{\begin{tabular}[t]{l}$\rho_1^{\vee}E$\end{tabular}}}}%
    \put(0.39575653,0.70592385){\color[rgb]{0,0,0}\makebox(0,0)[lt]{\lineheight{1.25}\smash{\begin{tabular}[t]{l}$\rho_3^{\vee}E$\end{tabular}}}}%
    \put(0.55685219,0.70592385){\color[rgb]{0,0,0}\makebox(0,0)[lt]{\lineheight{1.25}\smash{\begin{tabular}[t]{l}$\rho_{12}^{\vee}E$\end{tabular}}}}%
    \put(0.70549943,0.70592385){\color[rgb]{0,0,0}\makebox(0,0)[lt]{\lineheight{1.25}\smash{\begin{tabular}[t]{l}$\rho_{123}^{\vee}E$\end{tabular}}}}%
    \put(0.6673017,0.68357932){\color[rgb]{0,0,0}\makebox(0,0)[lt]{\lineheight{1.25}\smash{\begin{tabular}[t]{l}3\end{tabular}}}}%
    \put(0.20327055,0.68376312){\color[rgb]{0,0,0}\makebox(0,0)[lt]{\lineheight{1.25}\smash{\begin{tabular}[t]{l}1\end{tabular}}}}%
    \put(0.23089645,0.5739892){\color[rgb]{0,0,0}\makebox(0,0)[lt]{\lineheight{1.25}\smash{\begin{tabular}[t]{l}3\end{tabular}}}}%
    \put(0.47457266,0.58609565){\color[rgb]{0,0,0}\makebox(0,0)[lt]{\lineheight{1.25}\smash{\begin{tabular}[t]{l}2\end{tabular}}}}%
    \put(0.10061551,0.51468328){\color[rgb]{0,0,0}\makebox(0,0)[lt]{\lineheight{1.25}\smash{\begin{tabular}[t]{l}$\iota_0^{\vee}B$\end{tabular}}}}%
    \put(0.24245344,0.51468328){\color[rgb]{0,0,0}\makebox(0,0)[lt]{\lineheight{1.25}\smash{\begin{tabular}[t]{l}$\rho_1^{\vee}B$\end{tabular}}}}%
    \put(0.39575653,0.51468328){\color[rgb]{0,0,0}\makebox(0,0)[lt]{\lineheight{1.25}\smash{\begin{tabular}[t]{l}$\rho_3^{\vee}B$\end{tabular}}}}%
    \put(0.55685219,0.51468328){\color[rgb]{0,0,0}\makebox(0,0)[lt]{\lineheight{1.25}\smash{\begin{tabular}[t]{l}$\rho_{12}^{\vee}B$\end{tabular}}}}%
    \put(0.70549943,0.51468328){\color[rgb]{0,0,0}\makebox(0,0)[lt]{\lineheight{1.25}\smash{\begin{tabular}[t]{l}$\rho_{123}^{\vee}B$\end{tabular}}}}%
    \put(0.6673017,0.49233874){\color[rgb]{0,0,0}\makebox(0,0)[lt]{\lineheight{1.25}\smash{\begin{tabular}[t]{l}3\end{tabular}}}}%
    \put(0.20327055,0.49252255){\color[rgb]{0,0,0}\makebox(0,0)[lt]{\lineheight{1.25}\smash{\begin{tabular}[t]{l}1\end{tabular}}}}%
    \put(0.23089645,0.38274871){\color[rgb]{0,0,0}\makebox(0,0)[lt]{\lineheight{1.25}\smash{\begin{tabular}[t]{l}3\end{tabular}}}}%
    \put(0.47457266,0.39485531){\color[rgb]{0,0,0}\makebox(0,0)[lt]{\lineheight{1.25}\smash{\begin{tabular}[t]{l}2\end{tabular}}}}%
    \put(0.17571512,1.08682712){\color[rgb]{0,0,0}\makebox(0,0)[lt]{\lineheight{1.25}\smash{\begin{tabular}[t]{l}$\iota_1^{\vee}d$\end{tabular}}}}%
    \put(0.44441903,1.08682712){\color[rgb]{0,0,0}\makebox(0,0)[lt]{\lineheight{1.25}\smash{\begin{tabular}[t]{l}$\rho_2^{\vee}d$\end{tabular}}}}%
    \put(0.69269109,1.07603265){\color[rgb]{0,0,0}\makebox(0,0)[lt]{\lineheight{1.25}\smash{\begin{tabular}[t]{l}$\rho_{23}^{\vee}d$\end{tabular}}}}%
    \put(0.59696357,1.05330184){\color[rgb]{0,0,0}\makebox(0,0)[lt]{\lineheight{1.25}\smash{\begin{tabular}[t]{l}3\end{tabular}}}}%
    \put(0.34198141,1.05330184){\color[rgb]{0,0,0}\makebox(0,0)[lt]{\lineheight{1.25}\smash{\begin{tabular}[t]{l}2\end{tabular}}}}%
    \put(0.17571551,1.58405214){\color[rgb]{0,0,0}\makebox(0,0)[lt]{\lineheight{1.25}\smash{\begin{tabular}[t]{l}$\iota_1^{\vee}a$\end{tabular}}}}%
    \put(0.44441926,1.58405214){\color[rgb]{0,0,0}\makebox(0,0)[lt]{\lineheight{1.25}\smash{\begin{tabular}[t]{l}$\rho_2^{\vee}a$\end{tabular}}}}%
    \put(0.69269124,1.57325767){\color[rgb]{0,0,0}\makebox(0,0)[lt]{\lineheight{1.25}\smash{\begin{tabular}[t]{l}$\rho_{23}^{\vee}a$\end{tabular}}}}%
    \put(0.59696373,1.55052687){\color[rgb]{0,0,0}\makebox(0,0)[lt]{\lineheight{1.25}\smash{\begin{tabular}[t]{l}3\end{tabular}}}}%
    \put(0.34198172,1.55052687){\color[rgb]{0,0,0}\makebox(0,0)[lt]{\lineheight{1.25}\smash{\begin{tabular}[t]{l}2\end{tabular}}}}%
    \put(0.17571551,1.41831054){\color[rgb]{0,0,0}\makebox(0,0)[lt]{\lineheight{1.25}\smash{\begin{tabular}[t]{l}$\iota_1^{\vee}b$\end{tabular}}}}%
    \put(0.44441926,1.41831054){\color[rgb]{0,0,0}\makebox(0,0)[lt]{\lineheight{1.25}\smash{\begin{tabular}[t]{l}$\rho_2^{\vee}b$\end{tabular}}}}%
    \put(0.69269124,1.40751608){\color[rgb]{0,0,0}\makebox(0,0)[lt]{\lineheight{1.25}\smash{\begin{tabular}[t]{l}$\rho_{23}^{\vee}b$\end{tabular}}}}%
    \put(0.59696373,1.38478529){\color[rgb]{0,0,0}\makebox(0,0)[lt]{\lineheight{1.25}\smash{\begin{tabular}[t]{l}3\end{tabular}}}}%
    \put(0.34198172,1.38478529){\color[rgb]{0,0,0}\makebox(0,0)[lt]{\lineheight{1.25}\smash{\begin{tabular}[t]{l}2\end{tabular}}}}%
    \put(0.10061551,0.32344722){\color[rgb]{0,0,0}\makebox(0,0)[lt]{\lineheight{1.25}\smash{\begin{tabular}[t]{l}$\iota_0^{\vee}I$\end{tabular}}}}%
    \put(0.24245344,0.32344722){\color[rgb]{0,0,0}\makebox(0,0)[lt]{\lineheight{1.25}\smash{\begin{tabular}[t]{l}$\rho_1^{\vee}I$\end{tabular}}}}%
    \put(0.39575653,0.32344722){\color[rgb]{0,0,0}\makebox(0,0)[lt]{\lineheight{1.25}\smash{\begin{tabular}[t]{l}$\rho_3^{\vee}I$\end{tabular}}}}%
    \put(0.55685219,0.32344722){\color[rgb]{0,0,0}\makebox(0,0)[lt]{\lineheight{1.25}\smash{\begin{tabular}[t]{l}$\rho_{12}^{\vee}I$\end{tabular}}}}%
    \put(0.70549943,0.32344722){\color[rgb]{0,0,0}\makebox(0,0)[lt]{\lineheight{1.25}\smash{\begin{tabular}[t]{l}$\rho_{123}^{\vee}I$\end{tabular}}}}%
    \put(0.6673017,0.30110276){\color[rgb]{0,0,0}\makebox(0,0)[lt]{\lineheight{1.25}\smash{\begin{tabular}[t]{l}3\end{tabular}}}}%
    \put(0.20327055,0.30128656){\color[rgb]{0,0,0}\makebox(0,0)[lt]{\lineheight{1.25}\smash{\begin{tabular}[t]{l}1\end{tabular}}}}%
    \put(0.23089645,0.19151281){\color[rgb]{0,0,0}\makebox(0,0)[lt]{\lineheight{1.25}\smash{\begin{tabular}[t]{l}3\end{tabular}}}}%
    \put(0.47457266,0.20361925){\color[rgb]{0,0,0}\makebox(0,0)[lt]{\lineheight{1.25}\smash{\begin{tabular}[t]{l}2\end{tabular}}}}%
    \put(0.07362934,0.13220696){\color[rgb]{0,0,0}\makebox(0,0)[lt]{\lineheight{1.25}\smash{\begin{tabular}[t]{l}$\iota_0^{\vee}Z$\end{tabular}}}}%
    \put(0.24245344,0.13220696){\color[rgb]{0,0,0}\makebox(0,0)[lt]{\lineheight{1.25}\smash{\begin{tabular}[t]{l}$\rho_1^{\vee}Z$\end{tabular}}}}%
    \put(0.39575653,0.13220696){\color[rgb]{0,0,0}\makebox(0,0)[lt]{\lineheight{1.25}\smash{\begin{tabular}[t]{l}$\rho_3^{\vee}Z$\end{tabular}}}}%
    \put(0.55685219,0.13220696){\color[rgb]{0,0,0}\makebox(0,0)[lt]{\lineheight{1.25}\smash{\begin{tabular}[t]{l}$\rho_{12}^{\vee}Z$\end{tabular}}}}%
    \put(0.70549943,0.13220696){\color[rgb]{0,0,0}\makebox(0,0)[lt]{\lineheight{1.25}\smash{\begin{tabular}[t]{l}$\rho_{123}^{\vee}Z$\end{tabular}}}}%
    \put(0.6673017,0.10986242){\color[rgb]{0,0,0}\makebox(0,0)[lt]{\lineheight{1.25}\smash{\begin{tabular}[t]{l}3\end{tabular}}}}%
    \put(0.20327055,0.11004623){\color[rgb]{0,0,0}\makebox(0,0)[lt]{\lineheight{1.25}\smash{\begin{tabular}[t]{l}1\end{tabular}}}}%
    \put(0.23089645,0.00027239){\color[rgb]{0,0,0}\makebox(0,0)[lt]{\lineheight{1.25}\smash{\begin{tabular}[t]{l}3\end{tabular}}}}%
    \put(0.47457266,0.01237884){\color[rgb]{0,0,0}\makebox(0,0)[lt]{\lineheight{1.25}\smash{\begin{tabular}[t]{l}2\end{tabular}}}}%
  \end{picture}%
\endgroup%

%% file: figures/Imm_Fig8_zc.pdf_tex
\begingroup%
  \makeatletter%
  \providecommand\color[2][]{%
    \errmessage{(Inkscape) Color is used for the text in Inkscape, but the package 'color.sty' is not loaded}%
    \renewcommand\color[2][]{}%
  }%
  \providecommand\transparent[1]{%
    \errmessage{(Inkscape) Transparency is used (non-zero) for the text in Inkscape, but the package 'transparent.sty' is not loaded}%
    \renewcommand\transparent[1]{}%
  }%
  \providecommand\rotatebox[2]{#2}%
  \newcommand*\fsize{\dimexpr\f@size pt\relax}%
  \newcommand*\lineheight[1]{\fontsize{\fsize}{#1\fsize}\selectfont}%
  \ifx\svgwidth\undefined%
    \setlength{\unitlength}{200.21806612bp}%
    \ifx\svgscale\undefined%
      \relax%
    \else%
      \setlength{\unitlength}{\unitlength * \real{\svgscale}}%
    \fi%
  \else%
    \setlength{\unitlength}{\svgwidth}%
  \fi%
  \global\let\svgwidth\undefined%
  \global\let\svgscale\undefined%
  \makeatother%
  \begin{picture}(1,0.90448368)%
    \lineheight{1}%
    \setlength\tabcolsep{0pt}%
    \put(0,0){\includegraphics[width=\unitlength,page=1]{Imm_Fig8_zc.pdf}}%
    \put(0.81318078,0.89694698){\color[rgb]{1,0,0}\makebox(0,0)[lt]{\lineheight{1.25}\smash{\begin{tabular}[t]{l}$c_A(\xi_8^{\textit{tor}})$\end{tabular}}}}%
    \put(0.76857047,0.05525943){\color[rgb]{0.10196078,0.10196078,0.10196078}\makebox(0,0)[lt]{\lineheight{1.25}\smash{\begin{tabular}[t]{l}$c(\xi)$\end{tabular}}}}%
    \put(0.3622078,0.70852485){\color[rgb]{0.10196078,0.10196078,0.10196078}\makebox(0,0)[lt]{\lineheight{1.25}\smash{\begin{tabular}[t]{l}$x$\end{tabular}}}}%
    \put(0.91066696,0.1350019){\color[rgb]{0.10196078,0.10196078,0.10196078}\makebox(0,0)[lt]{\lineheight{1.25}\smash{\begin{tabular}[t]{l}$y$\end{tabular}}}}%
    \put(-0.00145954,0.11196413){\color[rgb]{0,0,1}\makebox(0,0)[lt]{\lineheight{1.25}\smash{\begin{tabular}[t]{l}$c_D(\xi^{\textit{in}}_8)$\end{tabular}}}}%
    \put(0.17125033,0.62922709){\color[rgb]{0.10196078,0.10196078,0.10196078}\makebox(0,0)[lt]{\lineheight{1.25}\smash{\begin{tabular}[t]{l}$z$\end{tabular}}}}%
  \end{picture}%
\endgroup%

%% file: figures/Imm_Fig8_zy.pdf_tex
\begingroup%
  \makeatletter%
  \providecommand\color[2][]{%
    \errmessage{(Inkscape) Color is used for the text in Inkscape, but the package 'color.sty' is not loaded}%
    \renewcommand\color[2][]{}%
  }%
  \providecommand\transparent[1]{%
    \errmessage{(Inkscape) Transparency is used (non-zero) for the text in Inkscape, but the package 'transparent.sty' is not loaded}%
    \renewcommand\transparent[1]{}%
  }%
  \providecommand\rotatebox[2]{#2}%
  \newcommand*\fsize{\dimexpr\f@size pt\relax}%
  \newcommand*\lineheight[1]{\fontsize{\fsize}{#1\fsize}\selectfont}%
  \ifx\svgwidth\undefined%
    \setlength{\unitlength}{200.21804449bp}%
    \ifx\svgscale\undefined%
      \relax%
    \else%
      \setlength{\unitlength}{\unitlength * \real{\svgscale}}%
    \fi%
  \else%
    \setlength{\unitlength}{\svgwidth}%
  \fi%
  \global\let\svgwidth\undefined%
  \global\let\svgscale\undefined%
  \makeatother%
  \begin{picture}(1,0.90551143)%
    \lineheight{1}%
    \setlength\tabcolsep{0pt}%
    \put(0,0){\includegraphics[width=\unitlength,page=1]{Imm_Fig8_zy.pdf}}%
    \put(0.81318087,0.89694707){\color[rgb]{1,0,0}\makebox(0,0)[lt]{\lineheight{1.25}\smash{\begin{tabular}[t]{l}$c_A(\xi_8^{\textit{tor}})$\end{tabular}}}}%
    \put(0.76857056,0.05525944){\color[rgb]{0.10196078,0.10196078,0.10196078}\makebox(0,0)[lt]{\lineheight{1.25}\smash{\begin{tabular}[t]{l}$c(\xi)$\end{tabular}}}}%
    \put(0.36220785,0.78197974){\color[rgb]{0.10196078,0.10196078,0.10196078}\makebox(0,0)[lt]{\lineheight{1.25}\smash{\begin{tabular}[t]{l}$x$\end{tabular}}}}%
    \put(0.91734483,0.13500192){\color[rgb]{0.10196078,0.10196078,0.10196078}\makebox(0,0)[lt]{\lineheight{1.25}\smash{\begin{tabular}[t]{l}$y$\end{tabular}}}}%
    \put(-0.00145951,0.11196414){\color[rgb]{0,0,1}\makebox(0,0)[lt]{\lineheight{1.25}\smash{\begin{tabular}[t]{l}$c_D(\xi^{\textit{in}}_8)$\end{tabular}}}}%
    \put(0.17125038,0.70268199){\color[rgb]{0.10196078,0.10196078,0.10196078}\makebox(0,0)[lt]{\lineheight{1.25}\smash{\begin{tabular}[t]{l}$z$\end{tabular}}}}%
  \end{picture}%
\endgroup%

%% file: figures/Imm_Fig8_xyc.pdf_tex
\begingroup%
  \makeatletter%
  \providecommand\color[2][]{%
    \errmessage{(Inkscape) Color is used for the text in Inkscape, but the package 'color.sty' is not loaded}%
    \renewcommand\color[2][]{}%
  }%
  \providecommand\transparent[1]{%
    \errmessage{(Inkscape) Transparency is used (non-zero) for the text in Inkscape, but the package 'transparent.sty' is not loaded}%
    \renewcommand\transparent[1]{}%
  }%
  \providecommand\rotatebox[2]{#2}%
  \newcommand*\fsize{\dimexpr\f@size pt\relax}%
  \newcommand*\lineheight[1]{\fontsize{\fsize}{#1\fsize}\selectfont}%
  \ifx\svgwidth\undefined%
    \setlength{\unitlength}{200.21806612bp}%
    \ifx\svgscale\undefined%
      \relax%
    \else%
      \setlength{\unitlength}{\unitlength * \real{\svgscale}}%
    \fi%
  \else%
    \setlength{\unitlength}{\svgwidth}%
  \fi%
  \global\let\svgwidth\undefined%
  \global\let\svgscale\undefined%
  \makeatother%
  \begin{picture}(1,0.9055115)%
    \lineheight{1}%
    \setlength\tabcolsep{0pt}%
    \put(0,0){\includegraphics[width=\unitlength,page=1]{Imm_Fig8_xyc.pdf}}%
    \put(0.81318078,0.89694708){\color[rgb]{1,0,0}\makebox(0,0)[lt]{\lineheight{1.25}\smash{\begin{tabular}[t]{l}$c_A(\xi_8^{\textit{tor}})$\end{tabular}}}}%
    \put(0.76857047,0.05525944){\color[rgb]{0.10196078,0.10196078,0.10196078}\makebox(0,0)[lt]{\lineheight{1.25}\smash{\begin{tabular}[t]{l}$c(\xi)$\end{tabular}}}}%
    \put(0.3622078,0.78197963){\color[rgb]{0.10196078,0.10196078,0.10196078}\makebox(0,0)[lt]{\lineheight{1.25}\smash{\begin{tabular}[t]{l}$x$\end{tabular}}}}%
    \put(0.91734462,0.13500191){\color[rgb]{0.10196078,0.10196078,0.10196078}\makebox(0,0)[lt]{\lineheight{1.25}\smash{\begin{tabular}[t]{l}$y$\end{tabular}}}}%
    \put(-0.00145954,0.11196404){\color[rgb]{0,0,1}\makebox(0,0)[lt]{\lineheight{1.25}\smash{\begin{tabular}[t]{l}$c_D(\xi^{\textit{in}}_8)$\end{tabular}}}}%
    \put(0.17125033,0.70268177){\color[rgb]{0.10196078,0.10196078,0.10196078}\makebox(0,0)[lt]{\lineheight{1.25}\smash{\begin{tabular}[t]{l}$z$\end{tabular}}}}%
  \end{picture}%
\endgroup%